\UseRawInputEncoding
\documentclass[12pt]{article}
\usepackage{amsmath,amssymb,amsthm, amsfonts,soul,verbatim,cite}
\usepackage{hyperref}
\usepackage{graphicx}
\usepackage{url}
\usepackage{dsfont} 
\usepackage{xcolor}
\definecolor{red}{rgb}{1,0,0}

\definecolor{blue}{rgb}{0,0,1}

\definecolor{green}{rgb}{0,.6,0}

\usepackage{float}
\usepackage{tikz}
\usetikzlibrary{topaths,calc}
\usepackage{enumerate}
\usepackage{subcaption}

\newtheorem{theorem}{Theorem}[section]
\newtheorem{corollary}[theorem]{Corollary}
\newtheorem{lemma}[theorem]{Lemma}
\newtheorem{proposition}[theorem]{Proposition}

\newtheorem{observation}[theorem]{Observation}

\newtheorem{construction}[theorem]{Construction}

\newtheorem{definition}[theorem]{Definition}
\newtheorem{claim}[theorem]{Claim}

\newcommand{\sat}{\mathrm{sat}}
\newcommand{\ex}{\operatorname{ex}_3}

\title{The Saturation Spectrum of Berge Stars}
\author{Neal Bushaw\footnote{Virginia Commonwealth University.  \texttt{nobushaw@vcu.edu}}, Sean English\footnote{University of North Carolina Wilmington.  \texttt{EnglishS@uncw.edu}}, Emily Heath\footnote{California State Polytechnic University Pomona. \texttt{eheath@cpp.edu}}, Daniel P. Johnston\footnote{Trinity College.   \texttt{daniel.johnston@trincoll.edu}}, Puck Rombach\footnote{University of Vermont.  \texttt{puck.rombach@uvm.edu}}}

\date{}

\begin{document}
	
	\maketitle
	
	\begin{abstract}
		The forbidden subgraph problem is among the oldest in extremal combinatorics -- how many edges can an $n$-vertex $F$-free graph have? The answer to this question is the well-studied extremal number of $F$. Observing that every extremal example must be maximally $F$-free, a natural minimization problem is also studied -- how few edges can an $n$-vertex maximal $F$-free graph have? This leads to the saturation number of $F$. Both of these problems are notoriously difficult to extend to $k$-uniform hypergraphs for any $k\ge 3$.
		
		Barefoot et al., in the case of forbidding triangles in graphs, asked a beautiful question -- which numbers of edges, between the saturation number and the extremal number, are actually realized by an $n$-vertex maximal $F$-free graph? Hence named the saturation spectrum of $F$, this has since been determined precisely for several classes of graphs through a large number of papers over the past two decades.
		
		In this paper, we extend the notion of the saturation spectrum to the hypergraph context. Given a graph $F$ and a hypergraph $G$ embedded on the same vertex set, we say $G$ is a {\bf{Berge-$F$}} if there exists a bijection $\phi:E(F)\to E(G)$ such that $e\subseteq \phi(e)$ for all $e\in E(F)$.  We completely determine the saturation spectrum for $3$-uniform Berge-$K_{1,\ell}$ for $1\leq \ell\leq 4$, and for $\ell=5$ when $5\mid n$. We also determine all but a constant number of values in the spectrum for $3$-uniform Berge-$K_{1,\ell}$ for all $\ell\geq 5$. We note that this is the first result determining the saturation spectrum for any non-trivial hypergraph.
	\end{abstract}
	
	\section{Introduction}
	
	A $k$-graph will refer to a $k$-uniform hypergraph; that is, a collection $V(G)$ of vertices, and a collection $E(G)\subseteq\binom{V(G)}{k}$. Given a $k$-graph $G$ and a forbidden $k$-graph $F$, we say $G$ is \textbf{$F$-saturated} if $G$ does not contain a copy of any $F$ as a subgraph, but for any $e\in E(\overline{G})$, $G+e$ does contain a copy of $F$.
	
	Perhaps the most widely studied topic in extremal graph theory is the \textbf{Tur\'an number} or 
	\textbf{extremal number}, denoted $\mathrm{ex}_k(n,F)$, which is the maximum number of edges among all $F$-saturated $k$-graphs on $n$ vertices. Extremal numbers were first studied in the context of graphs by Mantel~\cite{Mantel}, who proposed, as an exercise, the problem of showing that $\mathrm{ex}_2(n,K_3)=\lfloor n^2/4\rfloor$. This was extended by Tur\'an~\cite{Turan}, who determined $\mathrm{ex}_2(n,K_\ell)$ for all cliques.
	
	At the other end of the spectrum, the \textbf{saturation number}, denoted $\mathrm{sat}_k(n,F)$, is the minimum number of edges in an $n$-vertex $F$-saturated $k$-graph. Erd\H{o}s, Hajnal and Moon~\cite{EHM} determined $\mathrm{sat}_2(n,K_\ell)$ for all graph cliques. This was extended by Bollob\'as~\cite{Bollobas} for all $k$-graph cliques.
	
	In general, extremal numbers and saturation numbers can be very far apart. Indeed, for graphs, the celebrated Erd\H{o}s-Stone Theorem~\cite{ES} gives us that 
	\begin{equation}\label{equation ESS}
		\mathrm{ex}_2(n,F)=\frac{\chi(F)-2}{\chi(F)-1}\binom{n}{2}+o(n^2),
	\end{equation}
	where $\chi(F)$ is the chromatic number of the forbidden graph $F$. On the other hand, a well known result of K\'aszonyi and Tuza~\cite{KT} shows that for any fixed graph $F$, the saturation number is at most linear in $n$; that is, there is a constant $c_F$ so that
	\[
	\mathrm{sat}_2(n,F)\leq c_Fn.
	\]
	Thus, for non-bipartite graphs we always see a large gap between the saturation number and extremal number (if $F$ is bipartite,~\eqref{equation ESS} only gives us that $\mathrm{ex}_2(n,F)=o(n^2)$). While there is no direct analogue of the Erd\H{o}s-Stone theorem for hypergraphs, we still often find a large difference between saturation numbers and extremal numbers. In particular, a highly influential result of Pikhurko~\cite{Pikhurko} shows that for any $k$-graph $F$,
	\[
	\mathrm{sat}_k(n,F)=O(n^{k-1}).
	\]
	On the other hand, for example, for the complete $k$-graph on $\ell$ vertices $K_\ell^{(k)}$, we have $\mathrm{ex}_k(n,K_{\ell}^k)=\Theta(n^k)$~\cite{d}.
	
	The purpose of this paper is to explore the interval between the saturation and extremal numbers. The \textbf{edge spectrum} or \textbf{saturation spectrum}, denoted $\mathrm{ES}_k(n,\mathcal{F})$, is the collection of all integers $m$ such that there exists an $n$-vertex $m$-edge $\mathcal{F}$-saturated $k$-graph. This notion was perhaps first studied by Barefoot, Casey, Fisher, Fraughnaugh and Harary~\cite{Barefoot}, who considered maximal triangle-free graphs. Since then, the saturation spectrum for all graph cliques has been determined~\cite{Amin}, as well as many other graphs (see, e.g., \cite{SS2,SS3,SS4,SS5}). One result of particular interest to us is that for graph stars, $\mathrm{ES}_2(n,K_{1,\ell})$ is complete~\cite{SS2} - i.e. that 
	\[
	\mathrm{ES}_2(n,K_{1,\ell})=[\mathrm{sat}_2(n,K_{1,\ell}),\mathrm{ex}_2(n,K_{1,\ell})].
	\]
	This seems to be atypical, as all other graphs for which the saturation spectrum is known have some non-realizable values between the saturation number and extremal number. 
	
	In this work, we will study the saturation spectrum for a natural family of $3$-graphs. To our knowledge, this is the first work on saturation spectra of hypergraphs. 
	
	As a generalization of a Berge path and Berge cycle in a hypergraph~\cite{Berge}, Gerbner and Palmer~\cite{GP} proposed the following definition. Given a graph $F$ and a $k$-graph $H$, we will say that $H$ is a \textbf{Berge-$F$} if there exists a bijection $\phi:E(F)\to E(H)$ such that $F$ and $H$ can be embedded on the same vertex set in such a way that $e\subseteq \phi(e)$ for all $e\in E(F)$. In slight abuse of notation, we will sometimes write Berge-$F$ to denote the family of all $k$-graphs which are Berge-$F$. Both saturation numbers and extremal numbers of Berge graphs have recieved considerable attention lately. (See e.g.~\cite{EGGMS, EKNS, FL, G, GMP, GS}.)
	
	Here, we focus on the saturation spectrum for Berge stars. Gerbner, Methuku and Palmer determined the extremal number for Berge stars in many cases. That is, they determined precisely when the size of the host graph is divisible by the size of the star, and give very close bounds for non-divisible hosts.
	\begin{theorem}[Theorem 12 in \cite{GMP}]\label{theorem extremal number} If $\ell>k+1$, then
		\[
		\mathrm{ex}_k(n,\text{Berge-}K_{1,\ell})\leq \binom{\ell}k\frac{n}{\ell},
		\]
		and this is sharp whenever $\ell$ divides $n$. If $\ell\leq k+1$, then
		\[
		\mathrm{ex}_k(n,\text{Berge-}K_{1,\ell})\leq\left\lfloor \frac{(\ell-1)n}{k}\right\rfloor,
		\]
		and this is sharp whenever $n$ is large enough.
	\end{theorem}
	
	On the other hand, Austhof and English determined the saturation number of Berge stars. Here, a precise bound is given regardless of divisibility.
	\begin{theorem}[Theorem 4 in \cite{AE}]\label{theorem saturation number}
		For all $k\geq 3$, $\ell\in\mathbb{N}$ and $n$ large enough,
		\[
		\sat_k(n,\text{Berge-}K_{1,\ell})=\min_{a\in [n]\mid \binom{a-1}{k-1}\leq \ell-2} \left\lceil \frac{(\ell-1)(n-a)}{k}\right\rceil +\binom{a}k.
		\]
	\end{theorem}
	
	We note that while the extremal and saturation numbers for Berge stars both grow on the order of $n$, the asymptotics of the two functions are different when $k>\ell+1$. Therefore, there is still a gap between the extremal and saturation numbers. 
	
	Our first main result determines all but a constant number of values in the saturation spectrum for any size Berge star in a $3$-graph.
	
	\begin{theorem}\label{theorem main L>=5}
		Let $\ell\geq 5$. There exists a constant $c=c(\ell)$ such that for all $n$ large enough, and for all
		\[
		m\in \left[\sat_3(n,\text{Berge-}K_{1,\ell}),\ex(n,\text{Berge-}K_{1,\ell})-c\right]
		\]
		there exists an $n$-vertex $m$-edge Berge-$K_{1,\ell}$-saturated $3$-graph.
	\end{theorem}
	
	Theorem~\ref{theorem main L>=5} follows from Theorems~\ref{theorem lower range} and~\ref{theorem upper range}. Our second main result determines the $3$-uniform saturation spectrum for Berge-$K_{1,5}$ exactly for large $n$, as long as $5\mid n$. 
	
	\begin{theorem}\label{theorem main L=5}
		For sufficiently large $n$ with $5\mid n$, there exists an $n$-vertex $m$-edge Berge-$K_{1,5}$-saturated $3$-graph if and only if 
		\[
		m\in \left[\left\lceil\frac{4n}{3}\right\rceil,2n-4\right]\cup \{2n\}.
		\]
	\end{theorem}
	
	Theorem~\ref{theorem main L=5} follows from Theorems~\ref{theorem lower range} and~\ref{theorem upper range for L=5} and Proposition~\ref{proposition no saturated near the top}.
	
	In addition, for completeness, we determine the entire spectrum for $3$-uniform Berge-$K_{1,\ell}$ when $\ell\leq 4$. In this case, the saturation number and extremal number differ by at most $2$. This result and its short proof are given in Section \ref{section small stars}.
	
	\subsection{Definitions, notation and organization}
	
	Given a graph $F$ and hypergraph $G$ embedded on the same vertex set, we say that \textbf{$F$ witnesses a Berge-$F$ in $G$} if there exists an injection $\phi:E(F)\to E(G)$ such that $e\in \phi(e)$ for all $e\in E(F)$. Given a $3$-graph $G$ and a vertex $v\in V(G)$, the \textbf{Berge degree} of $v$ in $G$, denoted $d^B_G(v)$, is the largest integer $\ell$ such that a copy $S$ of $K_{1,\ell}$ can be embedded into $V(G)$ with the center of $S$ at $v$, such that $S$ witnesses a Berge-$K_{1,\ell}$. Given such an injection, we refer to the witness vertices of the edges as \textbf{core leaves}. We note that $G$ is Berge-$K_{1,\ell}$-free if and only if the Berge degree of every vertex is at most $\ell-1$.  When the host graph is clear from context, we frequently omit the subscript and write $d^B(v)$ for $d^B_G(v)$. Given a $k$-graph $G$ and an edge $e\in E(\overline{G})$, we say that $G+e$ contains a \textbf{new} Berge-$F$ if there exists a Berge-$F$ in $G+e$ that contains the edge $e$. 
	
	Given two $k$-graphs $G$ and $H$, the \textbf{disjoint union} $G\sqcup H$ is the $k$-graph formed where $V(G)$ and $V(H)$ are taken to be disjoint sets, and $V(G\sqcup H)=V(G)\cup V(H)$, while $E(G\sqcup H)=E(G)\cup E(H)$. The (not necessarily disjoint) \textbf{union} $G\cup H$ is defined similarly, where $V(G\cup H)=V(G)\cup V(H)$, while $E(G\cup H)=E(G)\cup E(H)$. The $2$-graph $K_\ell^-$ will be the graph formed from the clique $K_{\ell}$ by deleting a single edge. Given a $3$-graph $H$ we use $L_H(v)$ (or just $L(v)$) to denote the \textbf{$H$-link} of the vertex $v\in V(H)$. This is an auxiliary $2$-graph with $V(L_H(v))=N_H(v)$, and $xy\in E(L_H(v))$ if $xyv\in E(H)$.
	
	All asymptotics and Bachmann-Landau/`Big Oh' notation are with respect to $n$. All parameters are assumed to be constant with respect to $n$ unless specifically stated otherwise.  Throughout, we use standard notation and terminology wherever possible (see, e.g., \cite{Berge}).
	
	The rest of the paper is organized as follows. In Section~\ref{section small stars}, we show that the saturation spectrum for small stars is complete. While these results follow quickly, we include them as evidence that Berge-$K_{1,5}$ is in some sense the first interesting case. In Section~\ref{section tools}, we introduce some definitions and a randomized construction which will be useful for our main results. In Section~\ref{section lower range}, we show that the lower range of the saturation spectrum for Berge-$K_{1,\ell}$ is complete. In Section~\ref{section upper range L>=5}, we show that the upper range of the saturation spectrum for Berge-$K_{1,\ell}$ is nearly complete; that is, it contains the entire range from the saturation number to the extremal number aside from possibly constantly many values. Finally, in Section~\ref{section upper range L=5}, we determine the exact upper range of the spectrum of Berge-$K_{1,5}$, provided that $5\mid n$.

	\subsection{The edge spectrum for small stars}\label{section small stars}
	
	When $\ell\leq 4$, the extremal number and saturation number of Berge-$K_{1,\ell}$ are very close together. This is because the degree and the Berge degree of a vertex are the same if the degree is $4$ or less. This can be seen by e.g. considering all the ways that four or fewer $3$-edges can intersect in a single vertex. 
	
	\begin{proposition}
		Let $1\leq \ell\leq 4$. For $n$ large enough, the $3$-uniform saturation spectrum of Berge-$K_{1,\ell}$ is complete. In particular, for $n$ large enough, we have the following.
		\begin{enumerate}[(A)]
			\item $\displaystyle\sat_3(n,\text{Berge-}K_{1,1})=\ex(n,\text{Berge-}K_{1,1})=0$.\label{item sat 1}
			\item $\displaystyle\sat_3(n,\text{Berge-}K_{1,2})=\ex(n,\text{Berge-}K_{1,2})=\left\lfloor\frac{n}{3}\right\rfloor$.\label{item sat 2}
			\item $\sat_3(n,\text{Berge-}K_{1,3})=\ex(n,\text{Berge-}K_{1,3})$ if $(n \mod 3)=1$, and $\sat_3(n,\text{Berge-}K_{1,3})=\ex(n,\text{Berge-}K_{1,3})-1$ otherwise.\label{item sat 3}
			\item $\sat_3(n,\text{Berge-}K_{1,4})=n-2$, while $\ex(n,\text{Berge-}K_{1,4})=n$.\label{item sat 4}
		\end{enumerate}
	\end{proposition}
	
	\begin{proof}
		We have that \eqref{item sat 1} follows immediately since any edge forms a Berge-$K_{1,1}$. For~\eqref{item sat 2} through~\eqref{item sat 4}, it is easily verified that for $2\leq \ell\leq 4$, the minimum
		\[
		\min_{a\in [n]\mid \binom{a-1}{2}\leq \ell-2} \left\lceil \frac{(\ell-1)(n-a)}{3}\right\rceil +\binom{a}3
		\]
		always occurs at $a=2$. From this, along with Theorems~\ref{theorem saturation number} and~\ref{theorem extremal number}, parts~\eqref{item sat 2} and~\eqref{item sat 3} follow from considering cases based on the value of $n\textrm{ mod }3$, while~\eqref{item sat 4} follows immediately.
		From~\eqref{item sat 1},~\eqref{item sat 2} and~\eqref{item sat 3}, we can immediately conclude that the saturation spectrum is complete for $1\leq \ell\leq 3$. For $\ell=4$, we need an $n$-vertex saturated example on $n-1$ edges. Let $G$ be a $3$-graph which contains $n-2$ vertices of degree $3$ and two isolated vertices. (It is well known that $3$-regular $3$-graphs exist on all vertex sets that are large enough, see e.g. \cite{AE}.) Let $u,v\in V(G)$ be the two vertices of degree $0$, and let $xyz\in E(G)$. Then let $G'$ be the graph with $V(G')=V(G)$ and $E(G')=(E(G)\setminus \{xyz\})\cup \{xuv,yzu\}$. Then $E(G')$ has $n-1$ edges, no vertices of degree $4$ or more (and thus is Berge-$K_{1,4}$-free), and all vertices except for $u$ and $v$ have degree $3$. Therefore, adding any edge to $G'$ will create a Berge-$K_{1,4}$. 
	\end{proof}

	\section{Tools}\label{section tools}
	
	Here we introduce some basic tools and results which will be helpful in later sections. Given a ($2$-uniform) graph $G$, let the \textbf{tree component number} $\mathrm{tree}(G)$ denote the number of components of $G$ which are trees.
	
	\begin{lemma}\label{lemma bound on berge degree based on tree number}
		Let $H$ be a hypergraph and $v\in V(H)$. Then
		\[
		d^B(v)=|N(v)|-\mathrm{tree}(L(v)).
		\]
	\end{lemma}
	
	\begin{proof}
		Let $B$ be the incidence graph of $L(v)$. That is, $B$ is a bipartite graph with one part being $V(L(v))$ and the other $E(L(v))$; two vertices in $B$ are adjacent if they are incident in $L(v)$. Note that a matching in $B$ corresponds exactly to a choice of core leaves for a Berge star centered at $v$. Letting $\nu(B)$ denote the size of the largest matching in $B$, we then have
		\[d^B(v)=\nu(B).\]
		
		Let $C$ be a component of $L(v)$ and let $B_C$ be the incidence graph of $C$ (once again, this is the bipartite graph with parts $V(C)$ and $E(H)\cap\binom{C}{2}$, with adjacencies between a vertex and all edges containing it). Note that $\nu$ is additive over components, and therefore it will suffice to show that 
		\[
		\nu(B_C)=\begin{cases}
			|V(C)|,&\text{ if }|E(C)|\geq |V(C)|,\\
			|V(C)|-1,&\text{ if }|E(C)|=|V(C)|-1.
		\end{cases}
		\]
		In the case that $C$ is a tree (i.e. $|E(C)|=|V(C)|-1$), then $\nu(B_C)\leq |V(C)|-1$ since one part of $B_C$ has size $|V(C)|-1$, and one can greedily build a matching of size $|V(C)|-1$ by choosing any leaf of $C$, matching it with the single edge incident to it, deleting that leaf, and repeating.
		
		Now assume $|E(C)|\geq |V(C)|$. Note that $\nu(B_C)\leq |V(C)|$ and therefore we need only find a matching of the appropriate size. Since $C$ is not a tree, it contains some cycle $x_1,x_2,\dots,x_\ell$. Start with the matching which matches $x_i$ to $x_ix_{i+1}$ (with indices taken modulo $\ell$) for each $i\in [\ell]$. Note that this current matching has the property that if $X$ is the set of vertices in the matching and $Y$ is the set of edges, then $Y\subseteq \binom{X}{2}$. From here we will expand this matching iteratively, retaining this property.
		
		On each step of the iteration, we choose an unmatched neighbor $x$ of a currently matched vertex $y$, and match $x$ to the edge $xy$. Since $C$ is connected, this process can only terminate once every vertex in $V(C)$ is matched.
	\end{proof}

	\begin{definition}
		Given a $3$-graph $G$, we say a vertex $v\in V(G)$ is \textbf{aggressively Berge-$K_{1,\ell}$-saturated of type I} if $d^B(v)=\ell-1$ and the vertices in $L(v)$ that are not in a tree component form a clique. We say a vertex $v$ is \textbf{aggressively Berge-$K_{1,\ell}$-saturated of type II} if $v$ is not aggressively Berge-$K_{1,\ell}$-saturated of type I, $d^B(v)=\ell-1$, and for any two vertices $x,y\in N(v)$ such that neither $x$ nor $y$ is in a tree component of $L(v)$, either $xy\in E(L(v))$ or at least one of $x$ or $y$ is aggressively Berge-$K_{1,\ell}$-saturated of type I. We will simply say a vertex $v$ is \textbf{aggressively Berge-$K_{1,\ell}$-saturated} if it is of Type I or Type II. We will say a $3$-graph $G$ is \textbf{aggressively Berge-$K_{1,\ell}$-saturated} if for any Berge-$K_{1,\ell}$-saturated graph $H$, the disjoint union $G\sqcup H$ is Berge-$K_{1,\ell}$-saturated.
	\end{definition}

	When it is clear from context that we are considering Berge-$K_{1,\ell}$, we may just call vertices and graphs \textbf{aggressively saturated}. 
	
	\begin{lemma}\label{lemma saturated vertices}
		Let $G$ be a $3$-graph and let $e\in E(\overline{G})$ be a non-edge of $G$. If $e$ contains a vertex that is aggressively Berge-$K_{1,\ell}$-saturated, then $G+e$ contains a new Berge-$K_{1,\ell}$.
	\end{lemma}
	
	\begin{proof}
		We consider cases based on whether or not $e$ contains a vertex of Type I.
		
		\textbf{Case 1:} $e$ contains a vertex $v$ that is aggressively saturated of Type I, say $e=vxy$. If $x,y\in N_H(v)$, then they must be either in the same tree component of $L_H(v)$ or in different components of $L_H(v)$. In either case, $\mathrm{tree}(L_{H+e}(v))=\mathrm{tree}(L_H(v))-1$, while $|N_{H+e}(v)|=|N_{H}(v)|$. If $x\in N_H(v)$ and $y\not\in N_H(v)$, then $|N_{H+e}(v)|=|N_{H}(v)|+1$, while $\mathrm{tree}(L_{H+e}(v))=\mathrm{tree}(L_{H}(v))$. Finally if both $x,y\not\in N_H(v)$, then we have $|N_{H+e}(v)|=|N_{H}(v)|+2$, and $\mathrm{tree}(L_{H+e}(v))=\mathrm{tree}(L_H(v))+1$. In all cases, by Lemma~\ref{lemma bound on berge degree based on tree number} we have
		\[
		d^B_{H+e}(v)=|N_{H+e}(v)|-\mathrm{tree}(L_{H+e}(v))\geq |N_H(v)|-\mathrm{tree}(L_{H}(v))+1=d^B_H(v)+1=\ell.
		\]
		
		\textbf{Case 2:} $e$ contains a vertex $v$ that is aggressively saturated of Type II, say $e=vxy$, but no vertices of Type I. If both $x,y\in N(v)$, then note that since neither $x$ nor $y$ is Type I and $xy\not\in E(L_H(v))$, at least one of $x$ or $y$ must be in a tree component of $L(v)$. Regardless of if both $x$ and $y$ are in the same tree component or different components, we can see that  $\mathrm{tree}(L_{H+e}(v))(v)=\mathrm{tree}(L_{H}(v))-1$, while $|N_{H+e}(v)|=|N_{H}(v)|$. The cases where at least one of $x$ or $y$ is not in $N_H(v)$ are identical to Case 1. Therefore, again by Lemma~\ref{lemma bound on berge degree based on tree number}, we have
		\[
		d^B_{H+e}(v)=|N_{H+e}(v)|-\mathrm{tree}(L_{H+e}(v))\geq |N_H(v)|-\mathrm{tree}(L_{H}(v))+1=d^B_H(v)+1=\ell.
		\]
	\end{proof}
	
	\begin{corollary}\label{corollary if all vertices are agg sat}
		Let $G$ be a $3$-graph. If every vertex in $G$ is aggressively Berge-$K_{1,\ell}$-saturated, then $G$ is aggressively Berge-$K_{1,\ell}$-saturated.
	\end{corollary}
	
	\begin{proof}
		Let $G'$ be a Berge-$K_{1,\ell}$-saturated graph and let $H:=G\sqcup G'$. First note that $H$ is Berge-$K_{1,\ell}$-free since $G'$ is Berge-$K_{1,\ell}$-free and every vertex in $G$ has Berge degree $\ell-1$ meaning $G$ is also Berge-$K_{1,\ell}$-free. As there are no edges containing vertices of both $G$ and $G'$, then $H$ is clearly also Berge-$K_{1,\ell}$-free. Now, consider $e\in E(\overline{H})$. If $e$ contains any vertex from $V(G)$, then by Lemma~\ref{lemma saturated vertices} $H+e$ contains a Berge-$K_{1,\ell}$. If not, $e$ is contained in $V(G')$, but since $G'$ is saturated, $H+e$ still contains a Berge-$K_{1,\ell}$.
	\end{proof}
	
	\begin{corollary}\label{corollary if non agg sat induce clique}
		Let $G$ be a $3$-graph and $A\subseteq V(G)$ be the set of vertices of $G$ which are not aggressively Berge-$K_{1,\ell}$-saturated. If every vertex in $A$ has Berge degree at most $\ell-1$ and $A$ induces a clique in $G$, then $G$ is Berge-$K_{1,\ell}$-saturated.
	\end{corollary}
	
	\begin{proof}
		Every vertex in $G$ has Berge degree at most $\ell-1$, and thus $G$ is Berge-$K_{1,\ell}$-free. Consider $e\in E(\overline{G})$. Since $A$ induces a clique in $G$, we must have that $e$ contains a vertex outside $A$; in particular, it contains a vertex which is aggressively saturated. By Lemma~\ref{lemma saturated vertices}, $G+e$ then contains a Berge-$K_{1,\ell}$.
	\end{proof}
	
	Finally, we note that the disjoint union of two aggressively saturated graphs is aggressively saturated; this is immediate from the definition.
	
	\begin{observation}\label{observation union of aggresively saturated graphs}
		If $G$ and $G'$ are both aggressively Berge-$K_{1,\ell}$-saturated $3$-graphs, then $G\sqcup G'$ is aggressively Berge-$K_{1,\ell}$-saturated.
	\end{observation}
	
	\subsection{The configuration model}
	
	A hypergraph $H$ is called \textbf{linear} when each pair of edges intersects in at most one vertex. The construction used in \cite{AE} for the saturation number involves a linear hypergraph with a prescribed degree sequence. We will need a similar construction for our purposes, with a suitable modification to allow us to create saturated hypergraphs with specified properties.
	
	It is known that linear $k$-uniform $d$-regular hypergraphs exist -- such hypergraphs can be viewed as incidence structures called $(n,\frac{dn}{k},d,k)$-configurations. Such configurations are known to exist~\cite{Configs}, and thus so do the corresponding hypergraphs. In~\cite{AE} this was extended to \textbf{nearly $d$-regular} linear hypergraphs, in which every degree is either $d$ or $d-1$.  In this section, we prove a variant of this result where we prescribe the vertex degrees more precisely.
	
	We note that this is essentially a problem of finding a linear hypergraph with a specific degree sequence.  Characterizations of such graphical degree sequences are well known for graphs; the Erd\H{o}s-Gallai Theorem gives an easy to check condition for the existence of a graph with given degree sequence \cite{EG}; for more recent work see, e.g., \cite{RecentEG}.  However, no such result is known for $k$-uniform hypergraphs.  Recent work has explored such theorems in the case of non-uniform linear hypergraphs \cite{LinearEG}, but these results do not apply in our uniform case.  Specifically, we need the following somewhat technical result.
	
	\begin{lemma}\label{lemma random graph}
		Let $\ell\geq 5$ and $0\leq k\leq \binom{\ell}{3}$. Then there exists some $n_0=n_0(\ell)$ such that for all $n\geq n_0$, there exists a linear $3$-graph $G:=G(n,\ell,k)$ on $n$ vertices with the following properties.
		\begin{enumerate}[(I)]
			\item $G$ contains exactly $15k$ vertices of degree $\ell-5$.\label{property degree L-5}
			\item $G$ contains at most three vertices of degree $\ell-2$.\label{property degree L-2}
			\item All other vertices of $G$ have degree $\ell-1$.\label{property degree L-1}
			\item The $15k$ vertices of degree $\ell-5$ are pairwise non-adjacent.\label{property non-adjacent}
			\item There exists a pair of disjoint edges $e_1$ and $e_2$ of $G$ containing only vertices of degree $\ell-1$ such that no other edge $e'$ intersects both $e_1$ and $e_2$.\label{property disjoint edges}
		\end{enumerate}
	\end{lemma}

	To prove Lemma~\ref{lemma random graph}, we employ the probabilistic method. We describe below the \textbf{configuration model}, which gives a way to randomly generate a graph or hypergraph with a prescribed degree sequence. 
	
	A \textbf{pseudohypergraph} is a hypergraph that may contain loops and multiple edges. Formally, a pseudohypergraph is a pair $(V,E)$ where $V$ is a set of vertices, and $E$ is a multiset of \textbf{pseudoedges}, which itself is a multiset of vertices. Many standard notions which make sense for hypergraphs also translate naturally to pseudohypergraphs (such as degrees, uniformity, etc). We note that if all the pseudoedges of a pseudohypergraph have multiplicity $1$ and all the vertices in each pseudoedge have multiplicity $1$, we obtain a simple (non-pseudo) hypergraph.
	
	Let $n$ be an integer and $\mathbf{d}=(d_1,d_2,\dots,d_n)$ be a sequence of non-negative integers such that $3\mid \sum_{i=1}^n d_i$. We now describe a method to generate a random $3$-uniform pseudo-hypergraph on $n$ vertices with degree sequence $\mathbf{d}$. 
	
	Let $V_i=\{v_{i,j}\mid 1\leq j\leq d_i\}$ be a set of \textbf{configuration points} for each $1\leq i\leq n$, and let $S:=\bigcup_{i=1}^n V_i$. Let $M$ be a perfect $3$-matching on $S$, chosen uniformly at random among all such perfect matchings (note that the assumption  $3\mid \sum_{i=1}^n d_i$ guarantees that a perfect matching exists), which we will call a \textbf{configuration}. We then define $H$ to be the pseudo-$3$-graph on vertex set $V(H):=\{V_1,V_2\dots,V_n\}$, and for each $e=v_{i_1,j_1}v_{i_2,j_2}v_{i_3,j_3}\in M$, we add one copy of the multiset $V_{i_1}V_{i_2}V_{i_3}$ to $E(H)$. 
	
	We will denote the output of this process by $\mathbb{H}^{(3)}_{\mathbf{d}}$. We aim to show that $\mathbb{H}^{(3)}_{\mathbf{d}}$ has a positive probability of being a simple linear $3$-graph that satisfies all the properties of $G(n,\ell,k)$ from Lemma~\ref{lemma random graph}; this also shows that such a $3$-graph must exist.
	
	Let $\phi(x)$ denote the number of configurations on $x$ configuration points, and note that
	\begin{equation}\label{equation definition of phi x}
		\phi(x)=\frac{x!}{(3!)^{x/3}(x/3)!}.
	\end{equation}
	Given some $i\in [n]$, and $j_1,j_2\in [d_i]$, we will call the pair of configuration points $v_{i,j_1}$ and $v_{i,j_2}$ a \textbf{loop} if there is an edge of $M$ containing both $v_{i,j_1}$ and $v_{i,j_2}$. Similarly, given $i,i'\in [n]$, $j_1,j_2\in [d_{i}]$ and $j_1',j_2'\in [d_{i'}]$, we will call the quadruple $v_{i,j_1}$, $v_{i,j_2}$, $v_{i',j_1'}$, $v_{i',j_2'}$ an \textbf{overlap} if there exist edges $e_1$ and $e_2$ of $M$ such that $v_{i,j_1},v_{i',j_1'}\in e_1$ while $v_{i,j_2},v_{i',j_2'}\in e_2$. We note that if $M$ has no loops or overlaps, $\mathbb{H}^{(3)}_{\mathbf{d}}$ is a simple linear $3$-graph.
	
	The following probabilistic result allows us to bound the probability of loops and overlaps. Given integers $x$ and $s$, let $(x)_s$ denote the \textbf{falling factorial}, $(x)_s=\prod_{i=0}^{s-1}(x-i)$.
	
	\begin{theorem}[Theorem 6.10 in~\cite{JLR}]\label{theorem method of moments}
		Let $X_1,X_2,\dots,X_n,\dots$ and $Y_1,Y_2,\dots,Y_n,\dots$ be two sequences of random variables. If $\lambda,\mu\geq 0$ are real numbers such that as $n\to \infty$ we have
		\[
		\mathds{E}[(X_n)_{s_1}\cdot(Y_n)_{s_2}]\to \lambda^{s_1}\mu^{s_2},
		\]
		for all $s_1,s_2\in\mathbb{Z}_{\geq 0}$, then $X_n$ and $Y_n$ converge in distribution to independent Poisson random variables with mean $\lambda$ and $\mu$ respectively.
	\end{theorem}
	
	We can now begin making progress on the proof of Lemma~\ref{lemma random graph}. Let $\mathbf{d}(n,\ell,k)$ denote the degree sequence of length $n$ that has $15k$ vertices of degree $\ell-5$, at most $2$ vertices of degree $\ell-2$, all the rest degree $\ell-1$, such that $3\mid \sum_{i=1}^n d_i$.
	
	\begin{lemma}\label{lemma simple linear random graph}
		Let $\ell\geq 5$, $0\leq k\leq \binom{\ell}{3}$ and let $n$ be large. Then $\mathbb{H}^{(3)}_{\mathbf{d}(n,\ell,k)}$ is a simple linear $3$-graph with probability bounded away from $0$.
	\end{lemma}
	
	The proof of the above lemma is a variation on a standard argument about the configuration model and we will only provide a sketch of the proof. For complete details of how to prove this for the case $k=0$, see~\cite{AE}; the case for graphs appears in~\cite{JLR}.
	
	\begin{proof}[Proof sketch]
		Let $r$ be the remainder of $n(\ell-1)$ when divided by $3$, and note that we have exactly $N:=(n-15k)(\ell-1)+15k(\ell-4)-r$ configuration points. Recall the definition of $\phi(x)$ from~\eqref{equation definition of phi x}.
		
		We can then calculate the expected number of loops to be
		\[
		\sum_{i=1}^n\binom{d_i}{2}(N-2)\frac{\phi(N-3)}{\phi(N)}=\sum_{i=1}^n\binom{d_i}{2}\frac{2}{N-1}=n\binom{\ell-1}{2}\frac{2}{N-1}+o(1)=\ell-2+o(1).
		\]
		Similarly, we can calculate the expected number of overlaps to be
		\begin{align*}
			\sum_{i,j\in [n]} \binom{d_i}{2}\binom{d_j}{2}(N-4)(N-5)\frac{\phi(N-6)}{\phi(N)}&=\binom{n}{2}\binom{\ell-1}{2}^2\frac{4}{(N-2)(N-3)}+o(1)\\
			&=\frac{(\ell-2)^2}{2}+o(1).
		\end{align*}
		
		Based on this, we choose $\lambda :=\ell-2$ and $\mu:=\frac{(\ell-2)^2}{2}$. We omit the details, but with these, we may apply Theorem~\ref{theorem method of moments} to show that the number of loops and the number of overlaps converge to independent Poisson random variables with the given parameters. Since Poisson random variables takes value zero with positive probability (for $n$ large enough), the probability that $\mathbb{H}^{(3)}_{\mathbf{d}(n,\ell,k)}$ has no loops and no overlaps is bounded away from $0$.
	\end{proof}
	
	Our next lemma will show that an appropriate choice of $\mathbb{H}^{(3)}_{\mathbf{d}(n,\ell,k)}$ will also have Property~\ref{property non-adjacent} in Lemma~\ref{lemma random graph} with high probability.
	
	\begin{lemma}\label{lemma non-adjacency in our random graph}
		Let $\ell\geq 5$ and $0\leq k\leq \binom{\ell}{3}$. Then the degree $\ell-5$ vertices in $\mathbb{H}^{(3)}_{\mathbf{d}(n,\ell,k)}$ are pairwise non-adjacent with probability $1-o(1)$.
	\end{lemma}
	
	\begin{proof}
		Let $N$ be the number of configuration points for $\mathbb{H}^{(3)}_{\mathbf{d}(n,\ell,k)}$. Let $A$ be the collection of configuration points corresponding to vertices of degree $\ell-5$. Let $X$ be the random variable which counts the number of pairs $u,v\in A$ such that $u$ and $v$ end up in the same edge in our configuration $M$. Note that this counts not just adjacency among vertices of degree $\ell-5$, but also loops among these vertices. We have
		\[
		\mathds{E}[X]=\binom{A}{2}(N-2)\frac{\phi(N-3)}{\phi(N)}=\binom{A}2\frac{2}{N-1}=o(1)
		\]
		since $|A|$ is constant and $N\to \infty$ with $n$. Since $\mathds{E}[X]=o(1)$, we then must have $\Pr(X=0)=1-o(1)$.
	\end{proof}
	
	Lemmas~\ref{lemma simple linear random graph} and~\ref{lemma non-adjacency in our random graph} together will give us Properties~\eqref{property degree L-5} through~\eqref{property non-adjacent} with positive probability in Lemma~\ref{lemma random graph}. Property~\eqref{property disjoint edges} follows deterministically, as long as $n$ is large enough.
	
	\begin{claim}\label{claim random graph disjoint edges}
		Let $\ell\geq 5$, $0\leq k\leq \binom{\ell}{3}$ and $n$ large. Let $G$ be a simple $3$-graph with degree sequence $\mathbf{d}(n,\ell,k)$. Then $G$ contains a pair of disjoint edges $e_1$ and $e_2$ containing only vertices of degree $\ell-1$ such that no other edge $e'$ intersects both $e_1$ and $e_2$.
	\end{claim}
	
	\begin{proof}
		There are less than $(15k+2)(\ell-2)$ edges containing vertices that are not degree $\ell-1$, while $|E(G)|\geq \frac{(n-15k-2)(\ell-1)}{3}$. We can then choose an edge $e_1$ of $G$ containing only vertices of degree $\ell-1$.  Since $G$ has maximum vertex degree at most $\ell -1$, there are at most $3(\ell-2)$ edges that intersect $e_1$, and at most $2\cdot 3(\ell-2))^2$ edges that intersect these edges. As this number of intersecting edges depends only on $\ell$, we can ensure that we have a choice remaining for $e_2$ by letting $n$ grow sufficiently large.
	\end{proof}
	
	At last, we prove Lemma~\ref{lemma random graph}.
	
	\begin{proof}[Proof of Lemma~\ref{lemma random graph}]
		By Lemma~\ref{lemma simple linear random graph}, $\mathbb{H}^{(3)}_{\mathbf{d}(n,\ell,k)}$ is a simple linear $3$-graph that satisfies Properties~\eqref{property degree L-5} through~\eqref{property degree L-1} with probability bounded away from $0$. Then by Lemma~\ref{lemma non-adjacency in our random graph}, $\mathbb{H}^{(3)}_{\mathbf{d}(n,\ell,k)}$ satisfies Property~\eqref{property non-adjacent} with probability $1-o(1)$, and so as long as $n$ is large enough, $\mathbb{H}^{(3)}_{\mathbf{d}(n,\ell,k)}$ satisfies Properties~\eqref{property degree L-5} through~\eqref{property non-adjacent} with positive probability. Since such a graph can occur randomly from the configuration model, it must be that there exists some simple linear $3$-graph $G$ that has these properties. Finally, by Claim~\ref{claim random graph disjoint edges}, this $3$-graph $G$ also satisfies Property~\eqref{property disjoint edges}.
	\end{proof}
	
	\section{Lower range}\label{section lower range}
	
	Our goal in this section is to prove that when $m$ is near the saturation end of the spectrum, we can always find an $n$-vertex Berge-$K_{1,\ell}$-saturated $3$-graph on $m$ edges.
	\begin{theorem}\label{theorem lower range}
		Let $\ell\geq 5$ and $n$ large. For any integer 
		\[
		m\in \left[\mathrm{sat}_3(n,\text{Berge-}K_{1,\ell}),\frac{\ell(\ell-1)}{12}n\right],
		\]
		there exists an $n$-vertex Berge-$K_{1,\ell}$-saturated $3$-graph on $m$ edges.
	\end{theorem}
	
	We note that the choice of $\frac{\ell(\ell-1)}{12}n$ is somewhat arbitrary. The proof we use could be adapted for some $m>\frac{\ell(\ell-1)}{12}n$, but it cannot be extended all the way to $\ex(n,\text{Berge-}K_{1,\ell})$. Therefore, we simply choose a value that is convenient to cut off the lower range.
	
	\subsection{Constructions for the lower range}
	
	We provide a few constructions which will be necessary for the proof of Theorem~\ref{theorem lower range}, and prove some basic properties of the constructions. Our first construction is an aggressively saturated graph slightly sparser than $K^{(3)}_{\ell}$.
	
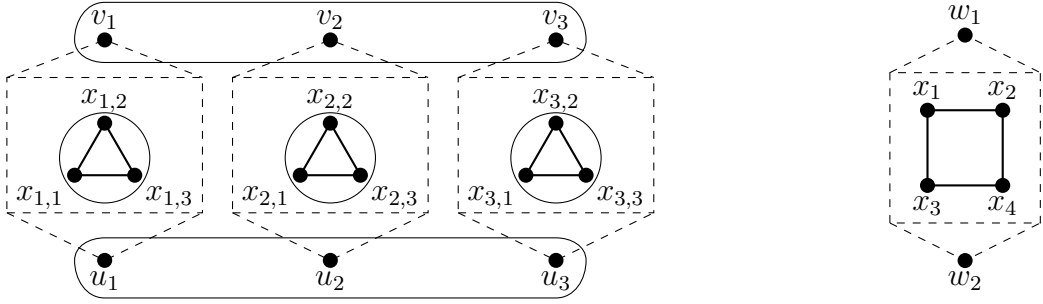
\begin{figure}
		\begin{subfigure}[t]{0.55\textwidth}
			\begin{center}
				\begin{tikzpicture}
					
					\draw[thick] (-0.4,0)--(0,.6928)--(.4,0)--(-0.4,0);
					\draw[thick] (2.6,0)--(3,.6928)--(3.4,0)--(2.6,0);
					\draw[thick] (5.6,0)--(6,.6928)--(6.4,0)--(5.6,0);
					
					\draw[dashed] (-1.3,-0.5)--(1.3,-0.5)--(1.3,1.3)--(-1.3,1.3)--(-1.3,-0.5);
					\draw[dashed] (1.7,-0.5)--(4.3,-0.5)--(4.3,1.3)--(1.7,1.3)--(1.7,-0.5);
					\draw[dashed] (4.7,-0.5)--(7.3,-0.5)--(7.3,1.3)--(4.7,1.3)--(4.7,-0.5);
					
					\draw[dashed] (-1.3,-0.5)--(0,-1.133)--(1.3,-0.5);
					\draw[dashed] (-1.3,1.3)--(0,1.8)--(1.3,1.3);
					
					\draw[dashed] (1.7,-0.5)--(3,-1.133)--(4.3,-0.5);
					\draw[dashed] (1.7,1.3)--(3,1.8)--(4.3,1.3);
					
					\draw[dashed] (4.7,-0.5)--(6,-1.133)--(7.3,-0.5);
					\draw[dashed] (4.7,1.3)--(6,1.8)--(7.3,1.3);

					\node (v1) at (0,1.8) {};
					\node (v2) at (3,1.8) {};
					\node (v3) at (6,1.8) {};
					\node (x11) at (-0.4,0) {};
					\node (x12) at (0,.6928) {};
					\node (x13) at (0.4,0) {};
					\node (x21) at (2.6,0) {};
					\node (x22) at (3,.6928) {};
					\node (x23) at (3.4,0) {};
					\node (x31) at (5.6,0) {};
					\node (x32) at (6,.6928) {};
					\node (x33) at (6.4,0) {};
					\node (u1) at (0,-1.133) {};
					\node (u2) at (3,-1.133) {};
					\node (u3) at (6,-1.133) {};
					\node (center1) at (0,.2309) {};
					\node (center2) at (3,.2309) {};
					\node (center3) at (6,.2309) {};
					
					\draw (center1) circle (0.6) node {};
					\draw (center2) circle (0.6) node {};
					\draw (center3) circle (0.6) node {};
					
					\fill (v1) circle (0.1) node [above] {$v_1$};
					\fill (v2) circle (0.1) node [above] {$v_2$};
					\fill (v3) circle (0.1) node [above] {$v_3$};
					\fill (x11) circle (0.1) node [below left] {$x_{1,1}$};
					\fill (x12) circle (0.1) node [above] {$x_{1,2}$};
					\fill (x13) circle (0.1) node [below right] {$x_{1,3}$};
					\fill (x21) circle (0.1) node [below left] {$x_{2,1}$};
					\fill (x22) circle (0.1) node [above] {$x_{2,2}$};
					\fill (x23) circle (0.1) node [below right] {$x_{2,3}$};
					\fill (x31) circle (0.1) node [below left] {$x_{3,1}$};
					\fill (x32) circle (0.1) node [above] {$x_{3,2}$};
					\fill (x33) circle (0.1) node [below right] {$x_{3,3}$};
					\fill (u1) circle (0.1) node [below] {$u_1$};
					\fill (u2) circle (0.1) node [below] {$u_2$};
					\fill (u3) circle (0.1) node [below] {$u_3$};
					
					\begin{scope}
						
						\draw ($(u1)-(0,0.5)$)
						to[out=0,in=180] ($(u3)-(0,0.5)$)
						to[out=0,in=270] ($(u3)+(0.4,0)$)
						to[out=90,in=0] ($(u3)+(0,0.3)$)
						to[out=180,in=0] ($(u1)+(0,0.3)$)
						to[out=180,in=90] ($(u1)-(0.4,0)$)
						to[out=270,in=180] ($(u1)-(0,0.5)$);
						
						\draw ($(v1)+(0,0.5)$)
						to[out=0,in=180] ($(v3)+(0,0.5)$)
						to[out=0,in=90] ($(v3)+(0.4,0)$)
						to[out=270,in=0] ($(v3)-(0,0.3)$)
						to[out=180,in=0] ($(v1)-(0,0.3)$)
						to[out=180,in=270] ($(v1)-(0.4,0)$)
						to[out=90,in=180] ($(v1)+(0,0.5)$);
					\end{scope}
					
				\end{tikzpicture}
			\end{center}
			\caption{The graph $L_5$ from Construction~\ref{construction lantern}, referred to as the lantern.}\label{figure L}
		\end{subfigure}~
		\begin{subfigure}[t]{0.45\textwidth}
			\begin{center}
				\begin{tikzpicture}
					
					\draw[thick] (0,0)--(0,1)--(1,1)--(1,0)--(0,0);
					
					\draw[dashed] (-0.5,-0.5)--(-0.5,1.5)--(1.5,1.5)--(1.5,-0.5)--(-0.5,-0.5);
					
					\draw[dashed] (-0.5,-0.5)--(0.5,-1)--(1.5,-0.5);
					\draw[dashed] (-0.5,1.5)--(0.5,2)--(1.5,1.5);

					\node (x1) at (0,1) {};
					\node (x2) at (1,1) {};
					\node (x3) at (0,0) {};
					\node (x4) at (1,0) {};
					\node (w1) at (0.5,2) {};
					\node (w2) at (0.5,-1) {};

					\fill (x1) circle (0.1) node [above] {$x_1$};
					\fill (x2) circle (0.1) node [above] {$x_2$};
					\fill (x3) circle (0.1) node [below] {$x_3$};
					\fill (x4) circle (0.1) node [below] {$x_4$};
					\fill (w1) circle (0.1) node [above] {$w_1$};
					\fill (w2) circle (0.1) node [below] {$w_2$};
					
				\end{tikzpicture}
			\end{center}
			\caption{The graph $S_5$ from Construction~\ref{construction S}, referred to as the sun graph.}\label{figure S}
		\end{subfigure}
		\caption{Diagrams of some of the constructions required for the proofs of Theorems~\ref{theorem lower range},~\ref{theorem upper range} and~\ref{theorem upper range for L=5}. A dashed rectangle represents that the $2$-graph inside the rectangle is part of the link of any vertex connected to the rectangle via dashed lines. All smooth curves represent $3$-edges in the graph.}\label{figure LandS}
	\end{figure}

	\begin{construction}\label{construction lantern}
		Let $\ell\geq 5$. The \textbf{$\ell$-lantern} $L_\ell$ is the graph on $3\ell$ vertices with
		\[
		V(L_\ell)=\{v_i,u_i,x_{i,j}\mid i\in [3],j\in [\ell-2]\}
		\]
		and edge set
		\[
		E(L_\ell)=\{v_1v_2v_3,u_1u_2u_3\}\cup \{v_ix_{i,j}x_{i,j'},u_ix_{i,j}x_{i,j'},x_{i,1}x_{i,2}x_{i,3}\mid i\in [3],j,j'\in [\ell-2], j\neq j'\}.
		\]
		We note that $|E(L_\ell)|=2+3\left(\binom{\ell-1}{3}+\binom{\ell-2}{2}\right)$.
		
		We will call the vertices of $L_\ell$ in $\{u_i,v_i\mid i\in [3]\}$ the \textbf{outer} lantern vertices, and the vertices in $\{x_{i,j}\mid i\in [3],j\in [\ell-2]\}$ the \textbf{inner} lantern vertices. See Figure~\ref{figure L} for a diagram of the $5$-lantern $L_5$.
	\end{construction}
	
	\begin{lemma}\label{lemma L is aggressively saturated}
		For all $\ell\geq 5$, $L_\ell$ is aggressively Berge-$K_{1,\ell}$-saturated.
	\end{lemma}
	
	\begin{proof}
		We will show that every vertex of $L_{\ell}$ is aggressively saturated. Then by Corollary~\ref{corollary if all vertices are agg sat}, $L_{\ell}$ is aggressively saturated. If $v$ is an outer lantern vertex, we can see that $L(v)\cong K_{\ell-2}\sqcup K_2$, and therefore $v$ is aggressively saturated of type I. If $v$ is an inner lantern vertex, then we find that $L(v)\cong K_{\ell-1}^-$, where the one missing edge in $L(v)$ contains two outer lantern vertices, and thus $v$ is aggressively saturated of type II. 
	\end{proof}
	
	We also observe easily that the clique $K_{\ell}^{(3)}$ is aggressively saturated.
	
	\begin{observation}\label{observation K is aggressively saturated}
		For all $\ell\geq 5$, $K_{\ell}^{(3)}$ is aggressively Berge-$K_{1,\ell}$-saturated.
	\end{observation}
	
	\begin{proof}
		The link of every vertex $v\in V(K_{\ell}^{(3)})$ is isomorphic to $K_{\ell-1}$. Therefore, each vertex is aggressively saturated of type I, and thus by Corollary~\ref{corollary if all vertices are agg sat}, $K_{\ell}^{(3)}$ is aggressively saturated.
	\end{proof}

	We now describe the main construction we will need for the lower range. The construction that achieves the sparsest Berge-$K_{1,\ell}$-saturated $3$-graph consists of a linear nearly-$(\ell-1)$-regular $3$-graph along with a small set of vertices forming a clique. To create constructions which are denser than this, we have three main tools: one tool corresponding to increasing the edge count of a Berge-$K_{1,\ell}$-saturated graph by large amounts, one tool corresponding to increasing by small amounts, and then one increasing by exactly one. 
	
	The tool which gives the largest changes in density comes from forming a disjoint union of a sparse graph with many copies of $K_\ell^{(3)}$, which are the densest Berge-$K_{1,\ell}$-free graphs. The second tool involves doing a local alteration to the sparse part of our construction -- small pockets of $15$ vertices will have slightly more edges than what we would have in a linear nearly-$(\ell-1)$-regular $3$-graph as each such pocket will correspond to having exactly three more edges than we would have had otherwise. Finally, we give a different local alteration to the sparse portion of our graph which can increase the number of edges by $1$ or $2$, allowing us to construct our desired hypergraph on any number of edges in the desired range.  The net effect of these alterations is summarized in the following explicit construction, and the alterations mentioned can be appropriately exploited by specifying the values $\ell,k,a$, and $i$.
	
	\begin{construction}\label{construction lower bound}
		Fix integers $\ell,k,a,i$ with $\ell\geq 5$, $0\leq k\leq \binom{\ell}{3}$, $3\leq a\leq \max\{3,\ell-3\}$, $0\leq i\leq 2$ and let $n$ be sufficiently large. We form a $3$-graph $W(n,\ell,k,a,i)$ as follows.
		
		Let $G:=G(n,\ell,k)$ from Lemma~\ref{lemma random graph}. Let $B$ be the set of $15k$ vertices of $G$ of degree $\ell-5$. Let $M\cong kL_5$ be a collection of $5$-lanterns with vertex set $V(M)=B$. Let $A$ be a copy of $K_a^{(3)}$ on the vertex set $V(A)=\{a_1,a_2,\dots,a_a\}$, disjoint from $V(G)$. Let us first form the auxiliary graph $G'$, where
		\[
		G'=\left(G\cup M\right)\sqcup A.
		\]
		Now, if $G$ contains one or two vertices of degree $\ell-2$, add an edge $e$ to $G'$ containing all of these vertices along with one or two vertices from $A$, and call the result $G''$. If $G$ does not contain vertices of degree $\ell-2$, then simply let $G'':=G'$.
		
		If $i=0$, let $W(n,\ell,k,a,i)=G''$. 
		
		If $i=1$, let $e_1$ and $e_2$ be edges of $G$ containing only vertices of degree (in $G$) $\ell-1$, such that no other edge of $G$ intersects both $e_1$ and $e_2$, say $e_1=x_1y_1z_2$ and $e_2=x_2y_2z_2$. Then let 
		$W(n,\ell,k,a,i)=G''-e_1-e_2+x_1x_2a_1+y_1y_2a_2+z_1z_2a_3$.
		
		If $i=2$, let $e_1$ be an edge in $G$ containing only vertices of degree (in $G$) $\ell-1$, say $e=xyz$. Then let $W(n,\ell,k,a,i)=G''-e_1+xa_1a_2+ya_2a_3+za_1a_3$.
	\end{construction}
	
	Let us now show that $W(n,\ell,k,a,i)$ is saturated.
	
	\begin{lemma}\label{lemma W is saturated}
		Fix integers $\ell\geq 5$, $0\leq k\leq \binom{\ell}{3}$, $3\leq a\leq \max\{3,\ell-3\}$ $0\leq i\leq 2$ and $n$ large. Then $W:=W(n,\ell,k,a,i)$ is Berge-$K_{1,\ell}$-saturated.
	\end{lemma}
	
	\begin{proof}
		Recalling the definitions of $G$, $M$, $A$, $e$, $e_1$, and $e_2$ from Construction~\ref{construction lower bound}, we will start by showing that every vertex in $V(G)=V(W)\setminus V(A)$ is aggressively saturated. We will consider the link graphs $L(v)$ for vertices $v\in V(G)$.
		
		If $v\in V(M)$, then we either have $L(v)\cong (\ell-4)K_2\sqcup K_3$ or $L(v)\cong (\ell-5)K_2\sqcup K_4^-$, depending on if $v$ is an outer vertex in a copy of $L_5$ or an inner vertex, respectively. In either case, using Lemma~\ref{lemma bound on berge degree based on tree number} we can calculate that $d^B(v)=\ell-1$. In the case where $v$ is an outer vertex of $L_5$, we see that $v$ is aggressively saturated of type I, and in the case where $v$ is an inner vertex, the one edge missing from the $K_4^-$-component of $L(v)$, call it $xy$, has the property that $x$ and $y$ are both outer vertices of $L_5$, which are aggressively saturated of type I, meaning $v$ is aggresively saturated of type II.
		
		If $v\in V(G)\cap e$, then $v$ was degree $\ell-2$ in $G$, and in particular we had $L_G(v)\cong (\ell-2)K_2$. The edge $e$ intersects $N_G(v)$ in at most one vertex (since at least one vertex in $e$ is in $A$). Therefore, depending on if $e\cap N_G(v)=\emptyset$ or not, we either have $L_W(v)\cong (\ell-1)K_2$ or $L_W(v)\cong (\ell-3)K_2\sqcup P_3$. In either case, using Lemma~\ref{lemma bound on berge degree based on tree number} we find that $d^B(v)=\ell-1$, and see that $v$ is aggressively saturated of type I.
		
		If $v\in e_1\cup e_2$, say, without loss of generality $v\in e_1$, then $L_G(v)\cong (\ell-1)K_2$. To form $W$, the edge $e_1$ was removed, and then a new edge containing $v$ and two vertices which were not in $N_G(v)$ was added. This yields $L_W(v)\cong (\ell-1)K_2$ and thus by Lemma~\ref{lemma bound on berge degree based on tree number}, $d^B(v)=\ell-1$ and so $v$ is aggressively saturated of type I.
		
		If $v$ is any vertex in $V(G)$ not covered above, then $v$ was degree $\ell-1$ in $G$, and therefore $L(v)\cong (\ell-1)K_2$. Hence by Lemma~\ref{lemma bound on berge degree based on tree number}, $d^B(v)=\ell-1$ and so $v$ is aggressively saturated of type I.
		
		We have exhausted all the possibilities and conclude that all vertices in $V(G)$ are aggressively saturated. In particular, no vertex in $V(G)$ can be the center of a Berge-$K_{1,\ell}$ since they all have Berge degree $\ell-1$. Furthermore, we claim that for any $v\in V(A)$ we have $d^B(v)<\ell$. Indeed, $d^B_A(v)=a-1\leq \ell-4$, and vertices from $A$ become incident with at most three more edges when forming $G''$ from $G'$, giving $v$ Berge degree at most $\ell-1$. Finally, $W$ is Berge-$K_{1,\ell}$-free. By Corollary~\ref{corollary if non agg sat induce clique} since the only non-aggressively saturated vertices in $W$ are in $V(A)$, and $V(A)$ induces a clique, we fiend that $W$ is Berge-$K_{1,\ell}$-saturated.
	\end{proof}
	
	We will need the following basic facts about $W(n,\ell,k,a,i)$ in order to precisely control the number of edges in our final construction.
	
	\begin{lemma}\label{lemma basic properties of W}
		Let $W:=W(n,\ell,k,a,i)$. The following hold.
		\begin{enumerate}[(a)]
			\item If $k=i=0$, then $|E(W)|=\left\lceil \frac{(\ell-1)(n-a)}{3}\right\rceil +\binom{a}3$.\label{statement saturated graph is in W}
			\item If $k<\binom{\ell}{3}$ and $W':=W(n,\ell,k+1,a,i)$, then $|E(W')|-|E(W)|=3$.\label{statement adding to k increases by 3}
			\item If $i\leq 1$ and $W':=W(n,\ell,k,a,i+1)$, then $|E(W')|-|E(W)|=1$.\label{statement adding to i increases by 1}
		\end{enumerate}
	\end{lemma}
	
	\begin{proof}
		Let us prove~\eqref{statement saturated graph is in W}. We note that $G:=G(n-a,\ell,0)$ is a linear nearly-$(\ell-1)$-regular $3$-graph and by a simple degree count, we can find that 
		\[
		E(G)= \left\lfloor \frac{(\ell-1)(n-a)}{3}\right\rfloor.
		\]
		If $\displaystyle \frac{(\ell-1)(n-a)}{3}$ is an integer, then this graph is $(\ell-1)$-regular, and so $W(n,\ell,0,a,0)$ has
		\[
		\displaystyle \frac{(\ell-1)(n-a)}{3}+\binom{a}{3}=\left\lceil \frac{(\ell-1)(n-a)}{3}\right\rceil +\binom{a}3
		\]
		edges. If $\displaystyle \frac{(\ell-1)(n-a)}{3}$ is not an integer, then $G$ contains a vertex of degree $\ell-2$. Therefore, in Construction~\ref{construction lower bound}, we add one extra edge from $V(G)$ to $V(A)$, along with the $\binom{a}{3}$ edges in $A$, giving us $\left\lfloor \frac{(\ell-1)(n-a)}{3}\right\rfloor+1+\binom{a}{3}=\left\lceil\frac{(\ell-1)(n-a)}{3}\right\rceil+\binom{a}{3}$ edges.
		
		Next we prove~\eqref{statement adding to k increases by 3}. Let $G:=G(n-a,\ell,k)$ and $G':=G(n-a,\ell,k+1)$ be the graphs used in the construction of $W$ and $W'$ respectively. Note that the degree sequences of $G$ and $G'$ differ only in that $G'$ has $15$ vertices of degree $\ell-5$, which are vertices of degree $\ell-1$ in $G$. This corresponds to the sum of degrees in $G'$ being $60$ less than the sum of degrees in $G$ and so $|E(G)|-|E(G')|=20$. When we then form $W$ and $W'$ from $G$ and $G'$, we add the same number of edges in every step, except we add a single extra copy of $L_5$ to $G'$. Since $|E(L_5)|=23$, this results in 
		\[
		|E(W')|-|E(W)|=|E(G')|-|E(G)|+23=3.
		\]
		
		Finally, we prove~\eqref{statement adding to i increases by 1}. The constructions for $W$ and $W'$ are identical except for the last step where we possibly remove edges $e_1$ and/or $e_2$, and add in a few new edges. When $i=0$, we do not change any edges; when $i=1$, we remove two edges, and add in three, resulting in one extra edge; when $i=2$, we remove one edge, and add in three, resulting in two extra edges. Thus, we indeed have $|E(W')|-|E(W)|=1$ when we increase $i$ by one.
	\end{proof}
	
	\subsection{Proof of \texorpdfstring{Theorem~\ref{theorem lower range}}{main lower range result}}
	
	Before we can prove our main result of this section, we need two basic facts about the saturation numbers for Berge-$K_{1,\ell}$. The first gives us some control over the value of $a$ in $W(n,\ell,k,a,i)$.
	
	\begin{claim}\label{claim a is not too big or too small}
		Let $\ell\geq 5$ and $n>\ell$. Then there exists 
		\[
		a^*\in \arg\min_{a\in [n]\mid \binom{a-1}{2}\leq \ell-2} \left\lceil \frac{(\ell-1)(n-a)}{3}\right\rceil +\binom{a}3
		\]
		such that
		\[
		a^*\in\begin{cases}
			\{3\}&\text{ if }\ell=5,\\
			[3,\ell-3]&\text{ if }\ell\geq 6.
		\end{cases}
		\]
	\end{claim}
	
	\begin{proof}
		Let $f(a):=\left\lceil \frac{(\ell-1)(n-a)}{3}\right\rceil +\binom{a}3$ and $M:=\arg\min_{a\in [n]\mid \binom{a-1}{2}\leq \ell-2} f(a)$. Note that
		\[
		\binom{a}2-\left\lceil \frac{\ell-1}{3}\right\rceil\leq f(a+1)-f(a)\leq \binom{a}2-\left\lfloor \frac{\ell-1}{3}\right\rfloor.
		\]
		For $a=1,2$ and $\ell\geq 5$, $\binom{a}2-\left\lfloor \frac{\ell-1}{3}\right\rfloor\leq 0$, and therefore there always exists some $a^*\in M$ with $a^*\geq 3$. Furthermore, when $a\geq 3$ and $\ell=5$, $\binom{a}2-\left\lceil \frac{\ell-1}{3}\right\rceil\geq 0$. Therefore, $3\in M$. 
		
		Now, we note that 
		\[
		\binom{a}{2}-\frac{\ell+1}{3}\leq \binom{a}2-\left\lceil \frac{\ell-1}{3}\right\rceil,
		\]
		and $\binom{a}{2}-\frac{\ell+1}{3}> 0$ for all $a>\frac{3+\sqrt{24\ell+33}}{6}$. In particular, when $\ell\geq 6$, we have that $\ell-3>\frac{3+\sqrt{24\ell+33}}{6}$, and so $a^*\leq \ell-3$ for all $a^*\in M$.
	\end{proof}
	
	Our second basic fact will help us bound one saturation number by another.
	
	\begin{observation}\label{observation monotone}
		Let $\ell$ be fixed. For all $n$ large enough, $\mathrm{sat}_3(n,\text{Berge-}K_{1,\ell})$ is monotone increasing in $n$.
	\end{observation}
	
	\begin{proof}
		This follows immediately from the expression given in Theorem~\ref{theorem saturation number}. Indeed, for any fixed $a$,
		\[
		\left\lceil \frac{(\ell-1)(n-\ell-a)}{3}\right\rceil +\binom{a}3
		\]
		is increasing in $n$. For $n$ large enough, since the possible values for $a$ are independent of $n$, increasing $n$ cannot decrease this minimum.
	\end{proof}

	We are now ready to prove Theorem~\ref{theorem lower range}.
	
	\begin{proof}[Proof of Theorem~\ref{theorem lower range}]
		Our construction will have the form $W\sqcup cK_\ell^{(3)}$, where $W=W(n-c\ell,\ell,k,a,i)$ for some suitably chosen $c,k,a,i$.
		Let $c$ be the largest integer such that
		
		\begin{equation}\label{equation definition of c}
			m-\binom{\ell}{3}c-\sat_3(n-c\ell,\text{Berge-}K_{1,\ell})\geq 0.
		\end{equation}
		
		First let us confirm that our choice of $c$ makes sense. We have that $c\geq 0$ since $m\geq\sat_3(n,\text{Berge-}K_{1,\ell})$. Furthermore, if $c^*>\frac{n}{4(\ell-2)}$, then
		\[
		m-\binom{\ell}{3}c^*-\sat_3(n-c^*\ell,\text{Berge-}K_{1,\ell})<\frac{\ell(\ell-1)}{12}n-\binom{\ell}{3}\frac{n}{4(\ell-2)}=0.
		\]
		Thus, $0\leq c\leq \frac{n}{4(\ell-2)}$, and $n-c\ell\geq n-\frac{n\ell}{4(\ell-2)}>n/2$.
		
		Using Claim~\ref{claim a is not too big or too small}, let 
		\[
		a^*\in \begin{cases}
			\{3\}&\text{ if }\ell=5,\\
			[3,\ell-3]&\text{ if }\ell\geq 6.
		\end{cases}
		\]
		be such that
		\[
		\sat_3(n-c\ell,\text{Berge-}K_{1,\ell})=\left\lceil \frac{(\ell-1)(n-c\ell-a^*)}{3}\right\rceil +\binom{a^*}3.
		\]
		Define
		\begin{equation}\label{equation defn of s}
			s:=m-c\binom{\ell}{3}-\sat_3(n-c\ell,\text{Berge-}K_{1,\ell}).
		\end{equation}
		We claim that $0\leq s < \binom{\ell}{3}$. Indeed, by~\eqref{equation definition of c}, we have that $s\geq 0$. Furthermore, by our choice of $c$, we have that
		\begin{equation}\label{equation upper bound on s}
			m-\binom{\ell}{3}(c+1)-\sat_3(n-(c+1)\ell,\text{Berge-}K_{1,\ell})<0.
		\end{equation}
		By Observation~\ref{observation monotone}, we have that
		\[
		\sat_3(n-c\ell,\text{Berge-}K_{1,\ell})\geq \sat_3(n-(c+1)\ell,\text{Berge-}K_{1,\ell}),
		\]
		so using this and the definition of $s$ with~\eqref{equation upper bound on s}, we have $s-\binom{\ell}{3}<0$.
		
		Now, let $k\in \mathbb{Z}_{\geq 0}$ and $i\in \{0,1,2\}$ be such that $s=3k+i$. Clearly $k<\binom{\ell}{3}$. Let $W:=W(n-c\ell,\ell,k,a^*,i)$. Using the fact that $K_\ell^{(3)}$ is aggressively saturated, we have from Observation~\ref{observation union of aggresively saturated graphs} that $cK_\ell^{(3)}$ is aggressively saturated. Furthermore by Lemma~\ref{lemma W is saturated}, $W$ is saturated. Thus, $H:=W\sqcup cK_\ell^{(3)}$ is saturated.
		
		Our final goal is to show that $E(H)=m$. From the definition of $s$ in~\eqref{equation defn of s}, it will suffice to show that $E(W)=\sat_3(n-c\ell,\text{Berge-}K_{1,\ell})+s$. Let $W'=W(n-c\ell,\ell,0,a^*,0)$. By Lemma~\ref{lemma basic properties of W}~\eqref{statement saturated graph is in W}, we have that
		\[
		|E(W')|=\left\lceil \frac{(\ell-1)(n-c\ell-a^*)}{3}\right\rceil +\binom{a^*}3=\sat_3(n-c\ell,\text{Berge-}K_{1,\ell}).
		\]
		Then, by repeated applications of Lemma~\ref{lemma basic properties of W}~\eqref{statement adding to k increases by 3} and~\eqref{statement adding to i increases by 1}, we have that $|E(W)|-|E(W')|=3k+i=s$, giving us that $|E(H)|=m$.
		
	\end{proof}

	\section{Approximate upper range for \texorpdfstring{$\ell\geq 5$}{L>=5}}\label{section upper range L>=5}
	
	Our main goal in this section is to prove the following.
	
	\begin{theorem}\label{theorem upper range}
		Let $\ell\geq 5$. There exists some constant $c>0$ and some $n_1$ such that for all $n\geq n_1$ and all
		\[
		m\in \left[\frac{\ell(\ell-1)}{12}n,\ex(n,\text{Berge-}K_{1,\ell})-c\right],
		\]
		there exists a Berge-$K_{1,\ell}$-saturated graph on $n$ vertices and $m$ edges.
	\end{theorem}
	
	\subsection{Constructions for the upper range (\texorpdfstring{$\ell\geq 5$}{L>=5})}
	
	We need to build up a few constructions that will allow us to control the number of edges in our final Berge-$K_{1,\ell}$-saturated graph. The first constuction is a small relatively sparse aggressively saturated $3$-graph.
	
	\begin{construction}\label{construction S}
		Let $S_{\ell}$ denote the $\textbf{sun graph}$ on $2\ell-4$ vertices with
		\[
		V(S_{\ell})=\{w_i,x_j\mid i\in [\ell-3],j\in [\ell-1]\}
		\]
		with
		\[
		E(S_{\ell})=\{w_ix_jx_{j+1}\mid i\in [\ell-3],j\in [\ell-1]\},
		\]
		where the indices of $x_j$ are taken modulo $\ell-1$. Note that $|E(S_{\ell})|=(\ell-1)(\ell-3)$. See Figure~\ref{figure S} for a diagram of the sun graph $S_5$.
	\end{construction}
	
	\begin{lemma}\label{lemma S is aggressively saturated}
		For all $\ell\geq 5$, $S_{\ell}$ from Construction~\ref{construction S} is aggressively Berge-$K_{1,\ell}$-saturated.
	\end{lemma}
	
	\begin{proof}
		For vertices $w_i$, we note that $L(w_i)$ is isomorphic to the $2$-graph cycle $C_{\ell-1}$. In particular from Lemma~\ref{lemma bound on berge degree based on tree number}, we can see that $d^B(w_i)=\ell-1$. Similarly, for vertices $x_j$, we can see that the link of $x_j$ is isomorphic to $K_{2,\ell-3}$. Again by Lemma~\ref{lemma bound on berge degree based on tree number}, we can see that $d^B(x_j)=\ell-1$. This gives us that $S_{\ell}$ is Berge-$K_{1,\ell}$-free, and furthermore, every vertex in $S_{\ell}$ has Berge degree $\ell-1$. 
		
		If we choose any edge $e\in E(\overline{S_{\ell}})$, then $e$ necessarily contains either two vertices from $\{w_i\mid i\in [\ell-3]\}$ or two non-consective vertices from $\{x_j\mid j\in [\ell-1]\}$. In either case, $e$ contains two vertices which are not adjacent in $S_\ell$. Since these vertices have Berge degree $\ell-1$ in $S_{\ell}$, they will have Berge degree $\ell$ in $S_{\ell}+e$. Thus $S_{\ell}$ is Berge-$K_{1,\ell}$-saturated.
		
		Now, since $S_{\ell}$ is Berge-$K_{1,\ell}$-saturated, and every vertex in $S_{\ell}$ has Berge degree $\ell-1$, $S_{\ell}$ is aggressively saturated, since for any saturated $G$, if we set $H:=S_{\ell}\sqcup G$ and choose any edge $e\in E(\overline{H})$, $e$ would necessarily be either contained in one of $V(S_{\ell})$ or $V(G)$ (creating a Berge-$K_{1,\ell}$ as both graphs are saturated), or would contain a vertex from $V(S_{\ell})$ and a vertex from $V(G)$, which would increase the Berge degree of the vertex from $S_{\ell}$, creating a Berge-$K_{1,\ell}$.
	\end{proof}

	We will use the following constructions as gadgets to control the number of edges in our final construction.
	
	\begin{construction}\label{construction H_1 and H_2}
		Let $\ell\geq 5$ and let $n$ be divisible by $3\ell$ and large enough that $G_1:=G(n,\ell,3)$ and $G_2:=G(n,\ell,1)$ from Lemma~\ref{lemma random graph} exist. We define two $3$-graphs, $H_1(n,\ell)$ and $H_2(n,\ell)$, as follows.
		
		Let $B_1$ and $B_2$ denote the sets of size $45$ and $15$ containing the vertices of degree $\ell-5$ in $G_1$ and $G_2$ respectively. Let $M\cong 3L_5$ be a $3$-graph with $V(M)=B_1$  and $K\cong 3K_5^{(3)}$ be a collection of three complete graphs with $V(K)=B_2$. Let $H_1(n,\ell):=G_1\cup M$ and let $H_2(n,\ell):=G_2\cup K$.
	\end{construction}
	
	The following lemma establishes a few basic facts about $H_1(n,\ell)$ and $H_2(n,\ell)$.
	
	\begin{lemma}\label{lemma properties of H_1 and H_2}
		Let $\ell\geq 5$ and let $n$ be divisible by $3\ell$ and large enough such that $G(n,\ell,1)$ and $G(n,\ell,3)$ both exist. Then the graphs $H_1(n,\ell)$ and $H_2(n,\ell)$ are aggressively Berge-$K_{1,\ell}$-saturated. Furthermore,
		\[
		|E(H_2(n,\ell))|-|E(H_1(n,\ell))|=1.
		\]
	\end{lemma}
	
	\begin{proof}
		Let $H_1:=H_1(n,\ell)$ and $H_2:=H_2(n,\ell)$. We first note that since $n$ is divisible by $3$, we have that for any $k\in \mathbb{Z}_{\geq 0}$
		\[
		\frac{(\ell-1)(n-15k)+(\ell-5)15k}{3}
		\]
		is an integer, and thus $G(n,\ell,k)$ does not have any vertices of degree $\ell-2$. 
		
		Let us now focus on showing $H_2$ is aggressively saturated. For any $v\in V(H_2)\setminus B_2$, we have that the link $L(v)\cong (\ell-1)K_2$, and thus $v$ is aggressively saturated of type I. Similarly, for any $v\in B_2$, we have $L(v)\cong(\ell-5)K_2\sqcup K_4$. Again, $v$ is aggressively saturated of type I. By Corollary~\ref{corollary if all vertices are agg sat}, $H_2$ is then aggressively saturated.
		
		Now let us focus on $H_1$. For any vertex $v$ in $V(H_1)\setminus B_1$, we have $L(v)\cong (\ell-1)K_2$, and therefore $v$ is aggressively saturated of type I. If $v\in B_1$, then $v$ is either an outer vertex or an inner vertex in one of the $5$-lanterns that were added to $H_1$. In the first case, we have $L(v)\cong (\ell-4)K_2\sqcup K_3$, and therefore $v$ is aggressively saturated of type I. In the second case, we have $L(v)\cong (\ell-5)K_2\sqcup K_4^-$, where the missing edge in the $K_4^-$ includes outer vertices of the lantern, which are aggressively saturated of type I, and thus $v$ is aggressively saturated of type II. Thus, again by Corollary~\ref{corollary if all vertices are agg sat}, $H_1$ is aggressively saturated.
		
		For the final claim, we can directly calculate that
		\[
		|E(H_1)|=|E(G_1)|+3|E(L_5)|=\frac{(\ell-1)(n-45)+(\ell-5)45}{3}+3(23)=\frac{(\ell-1)n}{3}+9,
		\]
		while
		\[
		|E(H_2)|=|E(G_2)|+3|E(K_5^{(3)})|=\frac{(\ell-1)(n-15)+(\ell-5)15}{3}+3(10)=\frac{(\ell-1)n}{3}+10.
		\]
		Thus, $|E(H_2)|-|E(H_1)|=1$.
	\end{proof}
	
	\begin{construction}\label{construction of upper bound graph}
		Let $\ell\geq 5$. We will define a $3$-graph $U(n,n_0,\ell,s,i)$ as follows. 
		
		Let $n_0$ be a fixed integer divisible by $3\ell$ and large enough such that $G(n_0,\ell,1)$ and $G(n_0,\ell,3)$ from Lemma~\ref{lemma random graph} exist. Let $H_1:=H_1(n_0,\ell)$ and $H_2:=H_2(n_0,\ell)$ be the graphs described in Construction~\ref{construction H_1 and H_2}. Let $n$ be large, and let $r$ be the remainder of $n$ when divided by $\ell$. For ease of notation, set
		\begin{equation}\label{equation definition of alpha}
			\alpha=\alpha(\ell):=(2\ell-4)\binom{\ell}{3}-\ell(\ell-1)(\ell-3).
		\end{equation}
		Let $0\leq i\leq \alpha$, and
		\[
		0\leq s<\frac{n-r-\alpha n_0}{(2\ell-4)\ell}
		\]
		be integers. Finally, for ease of notation, set
		\[
		\beta=\beta(n,n_0,\ell,s):=\frac{n-r-\alpha n_0-s\ell(2\ell-4)}{\ell}.
		\]
		Note that $\beta$ is an integer by our choice of $r$ and the fact that $\ell\mid n_0$, and that $\beta>0$ by our choice of $s$. With all this in hand, let
		\[
		U(n,n_0,\ell,s,i):=K^{(3)}_r\sqcup iH_2\sqcup(\alpha-i)H_1\sqcup (s\ell)S_\ell\sqcup \beta K_\ell^{(3)}.
		\]
	\end{construction}
	
	We note that the parameter $n_0$ that appears in $G(n,n_0,\ell,s,i)$ is independent of $n$, and thus will be treated as constant for asymptotic purposes.
	
	\begin{lemma}\label{lemma basic properties of U}
		Let $n$, $n_0$, $\ell$, $s$, $i$, $\alpha$ and $\beta$ be as defined in Construction~\ref{construction of upper bound graph}. The following hold.
		\begin{enumerate}[(i)]
			\item $U(n,n_0,\ell,s,i)$ is Berge-$K_{1,\ell}$-saturated.\label{statement U is saturated}
			\item $|V(U(n,n_0,\ell,s,i))|=n$.\label{statement U has size n}
			\item If $i<\alpha$, then $|E(U(n,n_0,\ell,s,i+1))|=|E(U(n,n_0,\ell,s,i))|+1$.\label{statement increasing i increases by 1 in U}
			\item If $s<\frac{n-r-\alpha n_0}{(2\ell-4)\ell}-1$, then $|E(U(n,n_0,\ell,s+1,i))|=|E(U(n,n_0,\ell,s,i))|-\alpha$.\label{statement increasing s increases by mu in U}
		\end{enumerate}
	\end{lemma}
	
	\begin{proof}
		By Lemma~\ref{lemma properties of H_1 and H_2}, we have that $H_1$ and $H_2$ are aggressively saturated. By Lemma~\ref{lemma S is aggressively saturated} and Observation~\ref{observation K is aggressively saturated}, we have that $S_{\ell}$ and $K_{\ell}^{(3)}$ are also aggressively saturated. Thus, $iH_2\sqcup(\alpha-i)H_1\sqcup (s\ell)S_\ell\sqcup \beta K_\ell^{(3)}$ is aggressively saturated by Observation~\ref{observation union of aggresively saturated graphs}. Finally since $K_r^{(3)}$ is  saturated, we must also have $U(n,n_0,\ell,s,i)$ is saturated, yielding ~\eqref{statement U is saturated}.
		
		Statement~\eqref{statement U has size n} follows in a straightforward manner from the definition of $\beta$ and the definition of $U(n,n_0,\ell,s,i)$; $\beta$ was chosen precisely so that $U(n,n_0,\ell,s,i)$ is order $n$.
		
		For~\eqref{statement increasing i increases by 1 in U}, we note that the only difference between $U(n,n_0,\ell,s,i+1)$ and $U(n,n_0,\ell,s,i)$ is that $U(n,n_0,\ell,s,i+1)$ has one copy of $H_2$ in place of a copy of $H_1$ in $U(n,n_0,\ell,s,i)$, and by Lemma~\ref{lemma properties of H_1 and H_2}, we have $|E(H_2)|-|E(H_1)|=1$.
		
		Similarly, for~\eqref{statement increasing s increases by mu in U}, we note that $\beta(n,n_0,\ell,s+1)=\beta(n,n_0,\ell,s)-(2\ell-4)$. Therefore, the difference between $U(n,n_0,\ell,s+1,i)$ and $U(n,n_0,\ell,s,i)$ is that $U(n,n_0,\ell,s+1,i)$ has $\ell$ copies of $S_\ell$ in place of $2\ell-4$ copies of $K_\ell^{(3)}$ in $U(n,n_0,\ell,s,i)$, Since $|E(S_\ell)|=(\ell-1)(\ell-3)$, we have that
		\begin{align*}
			|E(U(n,n_0,\ell,s+1,i))|&=|E(U(n,n_0,\ell,s,i))|-(2\ell-4)\binom{\ell}{3}+\ell(\ell-1)(\ell-3)\\
			&=|E(U(n,n_0,\ell,s,i))|-\alpha.
		\end{align*}
	\end{proof}

	\subsection{Proof of \texorpdfstring{Theorem~\ref{theorem upper range}}{main upper range result}}

	Before we prove the main result of this section, we need the following basic claim about rational functions.
	
	\begin{claim}\label{claim fraction L bound}
		For all $\ell\geq 5$,
		\[
		\frac{\ell-4}{4(\ell^2-7\ell+13)\ell}<\frac{1}{(2\ell-4)\ell}.
		\]
	\end{claim}
	
	\begin{proof}
		We have that
		\[
		4(\ell^2-7\ell+13)\ell-(\ell-4)(2\ell-4)\ell=2\ell(\ell^2-8\ell+18),
		\]
		which has only one real root at $\ell=0$, and is positive for all $\ell>0$. Thus, for all $\ell\geq 5$,
		\[
		4(\ell^2-7\ell+13)\ell-(\ell-4)(2\ell-4)\ell>0.
		\]
		Furthermore, $4(\ell^2-7\ell+13)\ell$ has only one real root, at $\ell=0$, and $(2\ell-4)\ell$ has roots $0$ and $2$, i.e. no roots when $\ell\geq 5$. Thus, rearranging the above inequality yields the desired result.
	\end{proof}
	
	\begin{proof}[Proof of Theorem~\ref{theorem upper range}]
		Choose $n_0$ divisible by $3\ell$ large enough so that $G(n_0,\ell,1)$ and $G(n_0,\ell,3)$ exist. Let $r$ be the remainder of $n$ when divided by $\ell$ and recall the definition of $\alpha$ from~\eqref{equation definition of alpha}. Let
		\[
		c:=\frac{\alpha n_0+\ell}{\ell}\binom{\ell}{3}.
		\]
		Note that by Theorem~\ref{theorem extremal number} and the fact that $r\leq \ell$,
		\begin{align*}
			|E(U(n,n_0,\ell,0,0))|&=\frac{n-r-\alpha n_0}{\ell}\binom{\ell}{3}+\binom{r}{3}+\alpha |E(H_1)|\\
			&\geq \ex(n,\text{Berge-}K_{1,\ell})-\frac{\alpha n_0+r}{\ell}\binom{\ell}{3}+\binom{r}{3}+\alpha|E(H_1)|\\
			&\geq \ex(n,\text{Berge-}K_{1,\ell})-c.
		\end{align*}
		With this in hand, let $m\in \left[\frac{\ell(\ell-1)}{12}n,\ex(n,\text{Berge-}K_{1,\ell})-c\right]$, and let
		\[
		m^*:=|E(U(n,n_0,\ell,0,0))|-m\geq 0.
		\]
		Then choose $a\geq 0$ and $0\leq b<\alpha$ so that $m^*=a\alpha-b$. We will provide an upper bound on $a$. We have
		\[
		m^*=|E(U(n,n_0,\ell,0,0))|-m\leq \binom{\ell}{3}\frac{n}{\ell}-\frac{\ell(\ell-1)}{12}n=\frac{\ell^2-5\ell+4}{12}n,
		\]
		so
		\begin{equation}\label{equation bound on a in upper range general proof}
			a\leq \frac{m^*}{\alpha}+1\leq\frac{\ell^2-5\ell+4}{12[(2\ell-4)\binom{\ell}{3}-\ell(\ell-1)(\ell-3)]}n+1=\frac{\ell-4}{4(\ell^2-7\ell+13)\ell}n+1.
		\end{equation}
		Now, since $\frac{\ell+\alpha n_0}{(2\ell-4)\ell}$ is constant, by Claim~\ref{claim fraction L bound}, as long as $n$ is large enough, we have 
		\begin{equation}\label{equation second bound on a in upper range general proof}
			\frac{\ell-4}{4(\ell^2-7\ell+13)\ell}n+1<\frac{n-\ell-\alpha n_0}{(2\ell-4)\ell}\leq \frac{n-r-\alpha n_0}{(2\ell-4)\ell}.
		\end{equation}
		Then combining~\eqref{equation bound on a in upper range general proof} and~\eqref{equation second bound on a in upper range general proof}, we have $a< \frac{n-r-\alpha n_0}{(2\ell-4)\ell}$.
		
		With this, we know that $U(n,n_0,\ell,a,b)$ exists, and through repeated applications of Lemma~\ref{lemma basic properties of U}, we have that $U(n,n_0,\ell,a,b)$ is an $n$-vertex graph that is Berge-$K_{1,\ell}$-saturated, and
		\[
		|E(U(n,n_0,\ell,a,b))|=|E(U(n,n_0,\ell,0,0))|-(a\alpha-b)=m.
		\]
	\end{proof}

	\section{The exact upper range for \texorpdfstring{$\ell=5$}{L=5} when \texorpdfstring{$5\mid n$}{5|n}}\label{section upper range L=5}
	
	Throughout this section we will assume that $n$ is large and that $5\mid n$. Since $5\mid n$, Theorem~\ref{theorem extremal number} gives us that $\ex(n,\text{Berge-}K_{1,5})=2n$.
	
	In this section, we will prove two results. First, for any $m\in [\frac{5}{3}n,2n-5]$, there exists a Berge-$K_{1,5}$-saturated $3$-graph on $n$ vertices and $m$ edges (see Theorem~\ref{theorem upper range for L=5}). Second, for any $m\in [2n-4,2n-1]$, there is no such $3$-graph (see Proposition~\ref{proposition no saturated near the top}). These results together with Theorem~\ref{theorem lower range} will be enough to completely determine the saturation spectrum in this case.
	
	\subsection{Excluded points near the top for \texorpdfstring{$\ell=5$}{L=5}}
	
	We will prove the following.
	\begin{proposition}\label{proposition no saturated near the top}
		Let $5\mid n$ and $i\in [4]$. Then there are no Berge-$K_{1,5}$-saturated graphs on $n$ vertices with $\ex(n,\text{Berge-}K_{1,5})-i$ edges.
	\end{proposition}
	
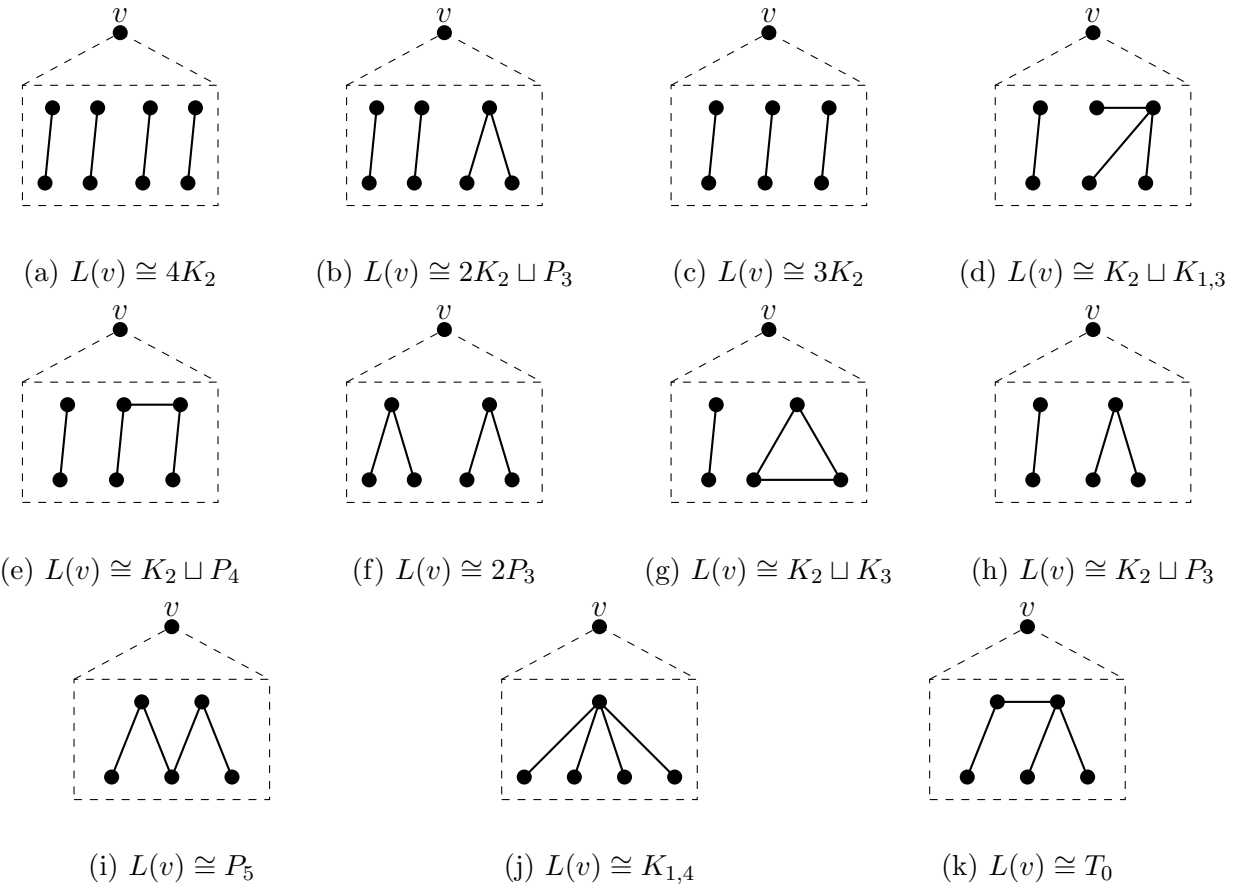
\begin{figure}
		\begin{subfigure}[t]{0.25\textwidth}
			\begin{center}
				\begin{tikzpicture}
					\draw[dashed] (-0.3,-0.3)--(2.3,-0.3)--(2.3,1.3)--(-0.3,1.3)--(-0.3,-0.3);
					\draw[dashed] (-0.3,1.3)--(1,2)--(2.3,1.3);
					
					\draw[thick] (0.1,1)--(0,0);
					\draw[thick] (.7,1)--(.6,0);
					\draw[thick] (1.4,1)--(1.3,0);
					\draw[thick] (2,1)--(1.9,0);
					
					\node (v) at (1,2) {};
					\node (1) at (0.1,1) {};
					\node (2) at (0,0) {};
					\node (3) at (.7,1) {};
					\node (4) at (.6,0) {};
					\node (5) at (1.4,1) {};
					\node (6) at (1.3,0) {};
					\node (7) at (2,1) {};
					\node (8) at (1.9,0) {};
					
					\fill (v) circle (0.1) node [above] {$v$};
					\fill (1) circle (0.1) node [below] {};
					\fill (2) circle (0.1) node [below] {};
					\fill (3) circle (0.1) node [below left] {};
					\fill (4) circle (0.1) node [below] {};
					\fill (5) circle (0.1) node [below] {};
					\fill (6) circle (0.1) node [below left] {};
					\fill (7) circle (0.1) node [below] {};
					\fill (8) circle (0.1) node [below] {};
				\end{tikzpicture}
			\end{center}
			\caption{$L(v)\cong 4K_2$}
		\end{subfigure}~
		\begin{subfigure}[t]{0.25\textwidth}
			\begin{center}
				\begin{tikzpicture}
					\draw[dashed] (-0.3,-0.3)--(2.3,-0.3)--(2.3,1.3)--(-0.3,1.3)--(-0.3,-0.3);
					\draw[dashed] (-0.3,1.3)--(1,2)--(2.3,1.3);
					
					\draw[thick] (0,0)--(0.1,1);
					\draw[thick] (0.6,0)--(0.7,1);
					\draw[thick] (1.9,0)--(1.6,1)--(1.3,0);
					
					\node (v) at (1,2) {};
					\node (1) at (0.1,1) {};
					\node (2) at (0,0) {};
					\node (3) at (.7,1) {};
					\node (4) at (.6,0) {};
					\node (5) at (1.6,1) {};
					\node (6) at (1.3,0) {};
					\node (8) at (1.9,0) {};
					
					\fill (v) circle (0.1) node [above] {$v$};
					\fill (1) circle (0.1) node [below] {};
					\fill (2) circle (0.1) node [below] {};
					\fill (3) circle (0.1) node [below left] {};
					\fill (4) circle (0.1) node [below] {};
					\fill (5) circle (0.1) node [below] {};
					\fill (6) circle (0.1) node [below left] {};
					\fill (8) circle (0.1) node [below] {};
				\end{tikzpicture}
			\end{center}
			\caption{$L(v)\cong 2K_2\sqcup P_3$}
		\end{subfigure}~
		\begin{subfigure}[t]{0.25\textwidth}
			\begin{center}
				\begin{tikzpicture}
					\draw[dashed] (-0.3,-0.3)--(2.3,-0.3)--(2.3,1.3)--(-0.3,1.3)--(-0.3,-0.3);
					\draw[dashed] (-0.3,1.3)--(1,2)--(2.3,1.3);
					
					\draw[thick] (0.3,1)--(0.2,0);
					\draw[thick] (1.05,1)--(0.95,0);
					\draw[thick] (1.8,1)--(1.7,0);
					
					\node (v) at (1,2) {};
					\node (1) at (0.3,1) {};
					\node (2) at (0.2,0) {};
					\node (5) at (1.05,1) {};
					\node (6) at (0.95,0) {};
					\node (7) at (1.8,1) {};
					\node (8) at (1.7,0) {};
					
					\fill (v) circle (0.1) node [above] {$v$};
					\fill (1) circle (0.1) node [below] {};
					\fill (2) circle (0.1) node [below] {};
					\fill (5) circle (0.1) node [below] {};
					\fill (6) circle (0.1) node [below left] {};
					\fill (7) circle (0.1) node [below] {};
					\fill (8) circle (0.1) node [below] {};
				\end{tikzpicture}
			\end{center}
			\caption{$L(v)\cong 3K_2$}
		\end{subfigure}~
		\begin{subfigure}[t]{0.25\textwidth}
			\begin{center}
				\begin{tikzpicture}
					\draw[dashed] (-0.3,-0.3)--(2.3,-0.3)--(2.3,1.3)--(-0.3,1.3)--(-0.3,-0.3);
					\draw[dashed] (-0.3,1.3)--(1,2)--(2.3,1.3);
					
					\draw[thick] (0.3,1)--(0.2,0);
					\draw[thick] (1.05,1)--(1.8,1)--(0.95,0);
					\draw[thick] (1.8,1)--(1.7,0);
					
					\node (v) at (1,2) {};
					\node (1) at (0.3,1) {};
					\node (2) at (0.2,0) {};
					\node (5) at (1.05,1) {};
					\node (6) at (0.95,0) {};
					\node (7) at (1.8,1) {};
					\node (8) at (1.7,0) {};
					
					\fill (v) circle (0.1) node [above] {$v$};
					\fill (1) circle (0.1) node [below] {};
					\fill (2) circle (0.1) node [below] {};
					\fill (5) circle (0.1) node [below] {};
					\fill (6) circle (0.1) node [below left] {};
					\fill (7) circle (0.1) node [below] {};
					\fill (8) circle (0.1) node [below] {};
				\end{tikzpicture}
			\end{center}
			\caption{$L(v)\cong K_2\sqcup K_{1,3}$}
		\end{subfigure}\\
		\begin{subfigure}[t]{0.25\textwidth}
			\begin{center}
				\begin{tikzpicture}
					\draw[dashed] (-0.3,-0.3)--(2.3,-0.3)--(2.3,1.3)--(-0.3,1.3)--(-0.3,-0.3);
					\draw[dashed] (-0.3,1.3)--(1,2)--(2.3,1.3);
					
					\draw[thick] (0.3,1)--(0.2,0);
					\draw[thick] (0.95,0)--(1.05,1)--(1.8,1);
					\draw[thick] (1.8,1)--(1.7,0);
					
					\node (v) at (1,2) {};
					\node (1) at (0.3,1) {};
					\node (2) at (0.2,0) {};
					\node (5) at (1.05,1) {};
					\node (6) at (0.95,0) {};
					\node (7) at (1.8,1) {};
					\node (8) at (1.7,0) {};
					
					\fill (v) circle (0.1) node [above] {$v$};
					\fill (1) circle (0.1) node [below] {};
					\fill (2) circle (0.1) node [below] {};
					\fill (5) circle (0.1) node [below] {};
					\fill (6) circle (0.1) node [below left] {};
					\fill (7) circle (0.1) node [below] {};
					\fill (8) circle (0.1) node [below] {};
				\end{tikzpicture}
			\end{center}
			\caption{$L(v)\cong K_2\sqcup P_4$}
		\end{subfigure}~
		\begin{subfigure}[t]{0.25\textwidth}
			\begin{center}
				\begin{tikzpicture}
					\draw[dashed] (-0.3,-0.3)--(2.3,-0.3)--(2.3,1.3)--(-0.3,1.3)--(-0.3,-0.3);
					\draw[dashed] (-0.3,1.3)--(1,2)--(2.3,1.3);
					
					\draw[thick] (0,0)--(0.3,1)--(0.6,0);
					\draw[thick] (1.9,0)--(1.6,1)--(1.3,0);
					
					\node (v) at (1,2) {};
					\node (1) at (0.3,1) {};
					\node (2) at (0,0) {};
					\node (4) at (.6,0) {};
					\node (5) at (1.6,1) {};
					\node (6) at (1.3,0) {};
					\node (8) at (1.9,0) {};
					
					\fill (v) circle (0.1) node [above] {$v$};
					\fill (1) circle (0.1) node [below] {};
					\fill (2) circle (0.1) node [below] {};
					\fill (4) circle (0.1) node [below] {};
					\fill (5) circle (0.1) node [below] {};
					\fill (6) circle (0.1) node [below left] {};
					\fill (8) circle (0.1) node [below] {};
				\end{tikzpicture}
			\end{center}
			\caption{$L(v)\cong 2P_3$}
		\end{subfigure}~
		\begin{subfigure}[t]{0.25\textwidth}
			\begin{center}
				\begin{tikzpicture}
					\draw[dashed] (-0.3,-0.3)--(2.3,-0.3)--(2.3,1.3)--(-0.3,1.3)--(-0.3,-0.3);
					\draw[dashed] (-0.3,1.3)--(1,2)--(2.3,1.3);
					
					\draw[thick] (0.2,0)--(0.3,1);
					\draw[thick] (.8,0)--(1.377,1)--(1.955,0)--(0.8,0);
					
					\node (v) at (1,2) {};
					\node (1) at (0.3,1) {};
					\node (2) at (0.2,0) {};
					\node (3) at (.8,0) {};
					\node (4) at (1.955,0) {};
					\node (5) at (1.377,1) {};
					
					\fill (v) circle (0.1) node [above] {$v$};
					\fill (1) circle (0.1) node [below] {};
					\fill (2) circle (0.1) node [below] {};
					\fill (3) circle (0.1) node [below left] {};
					\fill (4) circle (0.1) node [below] {};
					\fill (5) circle (0.1) node [below] {};
				\end{tikzpicture}
			\end{center}
			\caption{$L(v)\cong K_2\sqcup K_3$}
		\end{subfigure}~
		\begin{subfigure}[t]{0.25\textwidth}
			\begin{center}
				\begin{tikzpicture}
					\draw[dashed] (-0.3,-0.3)--(2.3,-0.3)--(2.3,1.3)--(-0.3,1.3)--(-0.3,-0.3);
					\draw[dashed] (-0.3,1.3)--(1,2)--(2.3,1.3);
					
					\draw[thick] (0.2,0)--(0.3,1);
					\draw[thick] (1,0)--(1.3,1)--(1.6,0);
					
					\node (v) at (1,2) {};
					\node (1) at (0.3,1) {};
					\node (2) at (0.2,0) {};
					\node (3) at (1,0) {};
					\node (4) at (1.6,0) {};
					\node (5) at (1.3,1) {};
					
					\fill (v) circle (0.1) node [above] {$v$};
					\fill (1) circle (0.1) node [below] {};
					\fill (2) circle (0.1) node [below] {};
					\fill (3) circle (0.1) node [below left] {};
					\fill (4) circle (0.1) node [below] {};
					\fill (5) circle (0.1) node [below] {};
				\end{tikzpicture}
			\end{center}
			\caption{$L(v)\cong K_2\sqcup P_3$}
		\end{subfigure}\\
		\begin{subfigure}[t]{0.3333\textwidth}
			\begin{center}
				\begin{tikzpicture}
					\draw[dashed] (-0.3,-0.3)--(2.3,-0.3)--(2.3,1.3)--(-0.3,1.3)--(-0.3,-0.3);
					\draw[dashed] (-0.3,1.3)--(1,2)--(2.3,1.3);
					
					\draw[thick] (0.2,0)--(0.6,1)--(1,0)--(1.4,1)--(1.8,0);
					
					\node (v) at (1,2) {};
					\node (1) at (1.4,1) {};
					\node (2) at (0.2,0) {};
					\node (3) at (.6,1) {};
					\node (4) at (1,0) {};
					\node (5) at (1.8,0) {};
					
					\fill (v) circle (0.1) node [above] {$v$};
					\fill (1) circle (0.1) node [below] {};
					\fill (2) circle (0.1) node [below] {};
					\fill (3) circle (0.1) node [below left] {};
					\fill (4) circle (0.1) node [below] {};
					\fill (5) circle (0.1) node [below] {};
				\end{tikzpicture}
			\end{center}
			\caption{$L(v)\cong P_5$}
		\end{subfigure}~
		\begin{subfigure}[t]{0.3333\textwidth}
			\begin{center}
				\begin{tikzpicture}
					\draw[dashed] (-0.3,-0.3)--(2.3,-0.3)--(2.3,1.3)--(-0.3,1.3)--(-0.3,-0.3);
					\draw[dashed] (-0.3,1.3)--(1,2)--(2.3,1.3);
					
					\draw[thick] (0,0)--(1,1)--(2,0);
					\draw[thick] (0.666,0)--(1,1)--(1.333,0);
					
					\node (v) at (1,2) {};
					\node (1) at (1,1) {};
					\node (2) at (0,0) {};
					\node (3) at (.666,0) {};
					\node (4) at (1.333,0) {};
					\node (5) at (2,0) {};
					
					\fill (v) circle (0.1) node [above] {$v$};
					\fill (1) circle (0.1) node [below] {};
					\fill (2) circle (0.1) node [below] {};
					\fill (3) circle (0.1) node [below left] {};
					\fill (4) circle (0.1) node [below] {};
					\fill (5) circle (0.1) node [below] {};
				\end{tikzpicture}
			\end{center}
			\caption{$L(v)\cong K_{1,4}$}
		\end{subfigure}~
		\begin{subfigure}[t]{0.3333\textwidth}
			\begin{center}
				\begin{tikzpicture}
					\draw[dashed] (-0.3,-0.3)--(2.3,-0.3)--(2.3,1.3)--(-0.3,1.3)--(-0.3,-0.3);
					\draw[dashed] (-0.3,1.3)--(1,2)--(2.3,1.3);
					
					\draw[thick] (1,0)--(1.4,1)--(1.8,0);
					\draw[thick] (0.2,0)--(0.6,1)--(1.4,1);
					
					\node (v) at (1,2) {};
					\node (1) at (1.4,1) {};
					\node (2) at (0.2,0) {};
					\node (3) at (.6,1) {};
					\node (4) at (1,0) {};
					\node (5) at (1.8,0) {};
					
					\fill (v) circle (0.1) node [above] {$v$};
					\fill (1) circle (0.1) node [below] {};
					\fill (2) circle (0.1) node [below] {};
					\fill (3) circle (0.1) node [below left] {};
					\fill (4) circle (0.1) node [below] {};
					\fill (5) circle (0.1) node [below] {};
				\end{tikzpicture}
			\end{center}
			\caption{$L(v)\cong T_0$}
		\end{subfigure}
		\caption{All allowable links of vertices $v$ with $|N(v)|\geq 5$ in a Berge-$K_{1,5}$-saturated graph.}\label{figure all allowable links}
	\end{figure}

	We will delay the proof of Proposition~\ref{proposition no saturated near the top} until the end of this section, after proving a number of ancillary results. Our first result characterizes the possible links we may see in a Berge-$K_{1,5}$-free graph. See Figure~\ref{figure all allowable links} for drawings of the links mentioned in the observation below. Let $T_0$ denote the unique tree of order $5$ with a degree $3$ vertex.
	
	\begin{observation}\label{observation allowable links}
		Let $G$ be a Berge-$K_{1,5}$-free graph and let $v\in V(G)$ be a non-isolated vertex. Then $2\leq |N(v)|\leq 8$ and the following hold.
		\begin{itemize}
			\item If $|N(v)|=8$, then $L(v)\cong 4K_2$.
			\item If $|N(v)|=7$, then $L(v)\cong 2K_2\sqcup P_3$.
			\item If $|N(v)|=6$, then $L(v)$ is isomorphic to a graph in $\{3K_2,K_2\sqcup K_{1,3},K_2\sqcup P_4,2P_3\}$.
			\item If $|N(v)|=5$, then $L(v)$ is isomorphic to a graph in $ \{K_2\sqcup K_3, K_2\sqcup P_3,P_5,K_{1,4},T_0\}$.
			\item $d(v)\leq 6$.
			\item If $d(v)=6$, then $L(v)\cong K_4$.
			\item If $d(v)=5$, then $L(v)\cong K_4^-$.
		\end{itemize}
	\end{observation}
	
	\begin{proof}
		Let $v$ be a non-isolated vertex. By Lemma~\ref{lemma bound on berge degree based on tree number}, we have that $|N(v)|-\mathrm{tree}(v)\leq 4$. Since every tree in $L(v)$ needs at least two vertices in $|N(v)|$, we need $|N(v)|\leq 8$. Furthermore, since $v$ is not isolated, $|N(v)|\geq 2$. Now, let us consider cases based on the neighborhood size of $v$. 
		
		\textbf{Case 1:} $|N(v)|=8$. Then $\mathrm{tree}(v)\geq 4$, and the only way we can have $4$ trees on $8$ vertices is $4K_2$.
		
		\textbf{Case 2:} $|N(v)|=7$. Then $\mathrm{tree}(V)\geq 3$, and the only way we can have three trees in $7$ vertices is $2K_2\sqcup P_3$.
		
		\textbf{Case 3:} $|N(v)|=6$. Then $\mathrm{tree}(V)\geq 2$. If we have two or more $K_2$-components, then we need three $K_2$-components, and therefore we get $3K_2$. If we have exactly one $K_2$-component, this leaves us with four vertices, which must be covered by a single tree (since isolated vertices cannot be in links), giving us two options, $K_2\sqcup K_{1,3}$ and $K_2\sqcup P_4$. Finally, if we have no $K_2$-components, we must have two trees of size $3$ each, giving us $2P_3$.
		
		\textbf{Case 4:} $|N(v)|=5$. Then $\mathrm{tree}(v)\geq 1$. If we have a $K_2$-component, then the second component can be any connected graph on $3$ vertices, giving us $K_2\sqcup K_3$ and $K_2\sqcup P_3$. If we do not have a $K_2$-component, then $L(v)$ must be connected and thus must be a tree. There are three trees of order $5$, namely, $P_5$, $K_{1,4}$ and $T_0$.
		
		Now, we note that if $|N(v)|>4$, one may quickly verify that none of the links discussed above have $5$ or more edges, if $|N(v)|<4$, then $d(v)\leq \binom{3}{2}=3$, and if $|N(v)|=4$, we have $d(v)\leq \binom{4}{2}=6$. Thus, if we have $d(v)\in \{5,6\}$, we must have $|N(v)|=4$, and it is readily verified that we have $L(v)\cong K_4$ if $d(v)=6$, or $L(v)\cong K_4^-$ if $d(v)=5$.
	\end{proof}

	The main idea of the proof of Proposition~\ref{proposition no saturated near the top} will be a simple degree counting argument. Recall that when $n$ is divisible by $5$, we have 
	\[
	\mathrm{ex}(n,\text{Berge-}K_{1,5})=2n,
	\]
	with the extremal construction being $n/5$ copies of $K^{(3)}_5$. By Observation~\ref{observation allowable links}, the maximum degree in a Berge-$K_{1,5}$-free graph is $6$, and indeed, in our extremal construction every vertex has degree exactly $6$. As such, to understand graphs that are missing only a few edges compared to this extremal example, we define the \emph{degree deficiency} of a vertex $v$ in a graph $G$ by
	\[
	\mathrm{ddf}_G(v):=6-d_G(v).
	\]
	We also define the degree deficiency of a set $S\subseteq V(G)$ to be
	\[
	\mathrm{ddf}_G(S):=\sum_{v\in S}\mathrm{ddf}_G(v).
	\]
	We often suppress the subscript $G$ if the host graph is clear from context. We also write $\mathrm{ddf}(G):=\mathrm{ddf}_G(V(G))$. We note that if $G$ is a graph on $n$ vertices (where $5|n$) and $\mathrm{ddf}(G)\geq 13$, then $|E(G)|\leq \mathrm{ex}(n,\text{Berge-}K_{1,5})-5$ by a simple degree count.

	We often find enough degree deficiency in the neighborhood of a single vertex. In order to bound the total degree deficiency, we have the following result which allows us to quickly bound the degree deficiency of a vertex based on its link.
	
	\begin{observation}\label{observation degrees}
		Let $G$ be a Berge-$K_{1,5}$-free $3$-graph. The following hold.
		\begin{itemize}
			\item For any vertices $u,v\in V(G)$, $\deg_{L(u)}(v)=\deg_{L(v)}(u)$.
			\item If $v\in V(G)$ and $L(v)$ contains a vertex of degree $1$, then $d(v)\leq 4$.
			\item If $v\in V(G)$ and $L(v)$ contains a vertex of degree $2$, then $d(v)\leq 5$.
			\item If $v\in V(G)$ and $L(v)$ contains a vertex of degree $4$, then $d(v)\leq 4$.
		\end{itemize}
	\end{observation}
	
	\begin{proof}
		Note that $\deg_{L(u)}(v)=|\{e\in E(G)\mid u,v\in e\}|=\deg_{L(v)}(u)$. Now, for $L(v)$ to have $5$ or more edges, we have to have $|N(v)|\geq 4$, and by Observation~\ref{observation allowable links}, we can see that no link with $|N(v)|\geq 5$ has $5$ or more edges (See Figure~\ref{figure all allowable links}), and therefore the only possible link with $6$ edges is $K_4$ and the only possible link with $5$ edges is $K_4^-$. Since $K_4$ only contains vertices of degree $3$, and $K_4^-$ only contains vertices of degree $2$ and $3$, the result follows.
	\end{proof}
	
	We can use Observation~\ref{observation degrees} to quickly get bounds on the degree deficiency of a graph containing a vertex with a given link. We summarize the bounds we get from analyzing each of the links featured in Observation~\ref{observation allowable links}.
	
	\begin{table}[h]
		\begin{center}
			\begin{tabular}{|c|c|c||c|c|c|}
				\hline
				$|N(v)|$&$L(v)$&$\mathrm{ddf}(N[v])\geq$&$|N(v)|$&$L(v)$&$\mathrm{ddf}(N[v])\geq$\\
				\hline
				$8$&$4K_2$&$18$&$5$&$K_2\sqcup K_3$&$9$\\
				$7$&$2K_2\sqcup P_3$&$15$&$5$&$K_2\sqcup P_3$&$12$\\
				$6$&$3K_2$&$15$&$5$&$P_5$&$9$\\
				$6$&$K_2\sqcup K_{1,3}$&$$14$$&$5$&$K_{1,4}$&$12$\\
				$6$&$K_2\sqcup P_4$&$12$&$5$&$T_0$&$9$\\
				$6$&$2P_3$&$12$&&&\\
				\hline
			\end{tabular}
			\caption{Lower bounds on the degree deficiency of a graph based on the link of a vertex $v$. See Figure~\ref{figure all allowable links} for drawings of the links in this table.}\label{table ddf bounds}
		\end{center}
	\end{table}
	
	\begin{observation}
		Let $G$ be a Berge-$K_{1,5}$-free graph. The bounds for $ddf(N[v])$ presented in Table~\ref{table ddf bounds} hold.
	\end{observation}
	
	\begin{proof}
		Given a vertex $v$, let $L_i$ denote the number of vertices of degree $i$ in $L(v)$. Then by Observation~\ref{observation degrees}, we have
		\begin{equation}\label{equation ddf in a link}
			\mathrm{ddf}(N[v])\geq ddf(v)+2L_1+L_2+2L_4=6-|E(L(v))|+2L_1+L_2+2L_4.
		\end{equation}
		Applying~\eqref{equation ddf in a link} to each link in the table yields the result. See Figure~\ref{figure all allowable links} for drawings of the links in question.
	\end{proof}
	
	Our goal now will be to show that every component which is not isomorphic to $K_5^{(3)}$ has too high of degree deficiency to exist in a graph with nearly $\mathrm{ex}(n,\text{Berge-}K_{1,5})$ edges. We first show that degree $6$ vertices cannot be adjacent in non-$K_5^{(3)}$ components.
	
	\begin{claim}\label{claim no adjacent degree 6}
		Let $G$ be a Berge-$K_{1,5}$-saturated graph. If $a$ and $a'$ are both degree $6$ vertices which are contained in some shared edge, then $a$ and $a'$ lie in a component which is isomorphic to $K_5^{(3)}$.
	\end{claim}
	
	\begin{proof}
		Since $a$ and $a'$ are both degree $6$, we have $L(a)\cong L(a')\cong K_4$, and since they are adjacent, we have $N(a)=\{a',b_1,b_2,b_3\}$ and $N(a')=\{a,b_1,b_2,b_3\}$. In $L(b_3)$ we find the edges $aa',ab_1,ab_2,a'b_1,a'b_2$, and therefore $K_4^-$ is a subgraph of $L(b_3)$. The only links which have this are $K_4$ and $K_4^-$. In particular, $N(b_3)=\{a,a',b_1,b_2\}$. By symmetry, this also holds for $b_1$ and $b_2$. Thus, the vertices $a,a',b_1,b_2 $ and $b_3$ form a component of order 5, which must then be isomorphic to $K_5^{(3)}$ in order to be Berge-$K_{1,5}^{(3)}$-saturated.
	\end{proof}
	
	We now consider possible neighborhood sizes for vertices, eliminating the possibility of large neighborhoods in Berge-$K_{1,5}$ saturated graphs.
	
	\begin{claim}\label{claim no vertices with large neighborhood}
		Let $H$ be a component in a Berge-$K_{1,5}$-saturated $3$-graph with $|V(H)|\geq 10$. If $H$ contains a vertex $v$ with $|N(v)|\geq 5$, then $\mathrm{ddf}(H)\geq 13$.
	\end{claim}
	
	\begin{proof}
		If $H$ contains a vertex $v$ with $N(v)\geq 7$, then by Table~\ref{table ddf bounds}, we are done. Thus, let us assume $H$ has no such vertex. Let $v\in V(H)$ have $|N(v)|\in \{5,6\}$. Let $R:=\{v\in V(H)\setminus N[v]\mid d(v)=6\}$. By Claim~\ref{claim no adjacent degree 6}, no two vertices in $R$ are adjacent.
		
		\textbf{Case 1:} $|N(v)|=6$. Then from Table~\ref{table ddf bounds}, we can see that $\mathrm{ddf}(N[v])\geq 12$. Then, we are done unless $R=V(H)\setminus N[v]$. We have $6|R|$ edges incident with vertices in $R$, and the other two endpoints of these edges must be contained in $N(v)$. Since $|R|\geq |V(H)|-|N[v]|\geq 3$, this gives us at least $18$ pairs (counted with multiplicity) of vertices from $N(v)$ which are in an edge with a vertex from $R$, however this implies that the average vertex from $N(v)$ is in $6$ such pairs, which would then imply that some vertices in $N(v)$ have degree strictly greater than $6$, a contradiction.
		
		\textbf{Case 2:} There are no vertices with neighborhood size $6$ or more, and $|N(v)|=5$. Then from Table~\ref{table ddf bounds}, we have $\mathrm{ddf}(N[v])\geq 9$. Let $V(H)=N[v]\cup R\cup Z$ be a partition where $R$ is the collection of degree $6$ vertices which are not in $N[v]$ (and hence $Z$ is the collection of vertices with degree less than $6$, which are not in $N[v]$). We are done unless $\mathrm{ddf}(Z)\leq 3$. In particular, we may assume $|Z|\leq 3$, and consequentially $|R|\geq 10-|N[v]|-|Z|\geq 1$. 
		
		Let $r\in R$. Since $L(r)\cong K_4$, we must have $|N(r)|=4$, and in particular, $N(r)\cap N(v)\neq\emptyset$ since $|Z|<|N(r)|$. Let $x\in N(r)\cap N(v)$, and let us consider $L(x)$. We know that $d_{L(x)}(r)=3$, and that $v\in L(x)$, while $rv\not\in L(x)$. This gives us that $|N(x)|>4$, in particular $|N(x)|=5$ due to the assumptions on this case. By inspecting all possible links, Observation~\ref{observation allowable links} shows that the only link with $5$ vertices and a vertex of degree (exactly) $3$ is $T$, the unique tree on $5$ vertices with a degree $3$ vertex. Let $y\in N(x)$ be the vertex that is degree $2$ in $L(x)$, and let $a_1,a_2$ be the remaining two vertices in $N(x)\setminus \{x,y,r\}$. Since $rv\not\in L(x)$, we must have that $v$ is the degree $1$ vertex adjacent to $y$ in $L(x)$, and so $a_1$ and $a_2$ are degree $1$ vertices adjacent to $r$. By Observation~\ref{observation degrees}, we also have $d_{L(v)}(x)=1$. Furthermore, we note that $y\in N(v)\cap N(r)$, and so by a symmetric argument, we can conclude that $d_{L(v)}(y)=1$. Thus, $xy$ form a $K_2$-component in $L(v)$. 
		
		Note that $N(v)\cap N(r)=\{x,y\}$ since if there was a third vertex $z\in N(v)\cap N(r)$, then $z$ would also need to be in a $K_2$-component of $L(v)$, but a link with $5$ vertices cannot contain two $K_2$-components. In particular, $a_1,a_2\not\in N[v]$, but $a_1$ and $a_2$ are degree one in $L(x)$. By Observation~\ref{observation degrees}, $\deg(a_1),\deg(a_2)\leq 4$. This gives us that
		\[
		\mathrm{ddf}(H)\geq \mathrm{ddf}(N[v])+\mathrm{ddf}(a_1)+\mathrm{ddf}(a_2)\geq 9+2+2=13.
		\]
	\end{proof}
	
	We now focus on components which only have vertices with neighborhood sizes of $4$ or less.
	
	\begin{claim}\label{claim small neighborhoods gives ddf 2}
		Let $H$ be a non-$K_5^{(3)}$ component of a Berge-$K_{1,5}$-saturated graph $G$. If $|N(v)|\leq 4$ for all $v\in V(H)$, then $\mathrm{ddf}(v)\geq 2$ for all $v\in V(H)$.
	\end{claim}
	
	\begin{proof}
		Recalling that $\mathrm{ddf}(v)=6-d(v)$, we assume to the contrary that there exists a vertex $v\in V(H)$ with $d(v)\geq 5$ (and hence $\mathrm{ddf}(v)<2$).
		
		\textbf{Case 1:} $d(v)=6$. Then by Observation~\ref{observation allowable links}, $L(v)\cong K_4$. Let $N(v)=\{a,b,c,d\}$. Since each vertex in $N(v)$ is connected to all others through an edge containing $v$ and no vertex has more than four neighbors, we must have $N[v]=N[a]=N[b]=N[c]=N[d]$. This implies that $V(H)=\{v,a,b,c,d\}$, and the only Berge-$K_{1,5}$-saturated graph on $5$ vertices is $K_5^{(3)}$.
		
		\textbf{Case 2:} $d(v)=5$. Then by Observation~\ref{observation allowable links}, $L(v)\cong K_4^-$. Let $N(v)=\{a,b,c,d\}$, where $ad$ is the missing edge in $L(v)$. Adding the edge $vad$ does not create a Berge-$K_{1,5}$ with center $v$. Therefore, the center of the new Berge-$K_{1,5}$ must be at either $a$ or $d$, assume without loss of generality that the center is at $a$. Note that based off $L(v)$, $N_G(b)=\{v,a,c,d\}$ while $N_G(c)=\{v,a,b,d\}$, and thus $b,c\in N(a)$. Consider $L(a)$. We must have $|N(a)|=4$, since if $N(a)$ was any smaller, adding the edge $vad$ would not create a Berge-$K_{1,5}$ at $a$. Let $x$ be the fourth vertex in $N(a)$ (we already know that $v,b$ and $c$ are all in $N(a)$). Note that in $L(a)$, $x$ must be adjacent to some vertex in $N(a)\setminus\{x\}=\{v,b,c\}$. Then, we must have $x=d$. However, then adding the edge $vad$ does not create a Berge-$K_{1,5}$ centered at $a$, a contradiction.
	\end{proof}

	We are now ready to prove the main result of this section.
	
	\begin{proof}[Proof of Proposition~\ref{proposition no saturated near the top}]
		Let $G\not\cong \frac{n}{5}K_5^{(3)}$ be a Berge-$K_{1,5}$-saturated $3$-graph. Note that $G$ must have at least one component which is not isomorphic to $K_5^{(3)}$. We will show that $\mathrm{ddf}(G)\geq 13$, which will imply that
		\[
		|E(G)|\leq \ex(n,\text{Berge-}K_{1,5})-5.
		\]
		Let $H_1,H_2,\dots,H_k$ be all the non-$K_5^{(3)}$ components in $G$, and let $n_0:=\sum_{i=1}^k |V(H_i)|$. Since both $n$ and $|V(K_5^{(3)})|$ are divisible by five, we must have that $5|n_0$, and furthermore since $G$ is saturated, we could not have $n_0=5$; hence $n_0\geq 10$. We now break into cases depending on the number of non-$K_5^{(3)}$ components present.
		
		\textbf{Case 1:} $k=1$. Then we have a non-$K_5^{(3)}$ component, $H_1$ with $10$ or more vertices, and therefore by Claim~\ref{claim no vertices with large neighborhood}, $H_1$ does not have any vertices with neighborhood size $5$ or more. By Claim~\ref{claim small neighborhoods gives ddf 2}, this gives us that $\mathrm{ddf}(H_1)\geq 2|V(H_1)|\geq 20$.
		
		\textbf{Case 2:} $k\geq 2$. If no vertex in $\bigcup_{i=1}^k V(H_i)$ has neighborhood size greater than $4$, then by Claim~\ref{claim small neighborhoods gives ddf 2}, we have $\mathrm{ddf}(G)\geq 2n_0\geq 20$. If we have two vertices, $x$ and $y$, from different components with $|N(x)|,|N(y)|\geq 5$, then by Table~\ref{table ddf bounds}, we have $\mathrm{ddf}(N[x]),\mathrm{ddf}(N[y])\geq 9$, and therefore $\mathrm{ddf}(G)\geq \mathrm{ddf}(N[x])+\mathrm{ddf}(N[y])\geq 18$. Thus, the only remaining subcase is when there is a single component $H_1$ with a vertex, say $x\in V(H_1)$ with $|N(x)|\geq 5$, and all other $K_5^{(3)}$ components having neighborhoods of size at most $4$. If $|V(H_2)|=1$, then $H_2$ consists of a vertex of degree $0$, and therefore
		\[
		\mathrm{ddf}(G)\geq \mathrm{ddf}(N[x])+\mathrm{ddf}(H_2)\geq 9+6=15,
		\]
		and if $|V(H_2)|\geq 2$, then by Claim~\ref{claim small neighborhoods gives ddf 2}, we have
		\[
		\mathrm{ddf}(G)\geq \mathrm{ddf}(N[x])+\mathrm{ddf}(H_2)\geq 9+2|V(H_2)|\geq 13.
		\]
	\end{proof}
	
	\subsection{Constructions for the upper range (\texorpdfstring{$\ell=5$}{L=5})}
	
	The reason the proof in the case $\ell=5$ is tractile stems from the fact that the $5$-lantern $L_5$ is aggressively saturated and not much sparser than $K_5^{(3)}$. In particular, in the extremal example $\frac{n}{5}K_5^{(3)}$, if we replace $3$ copies of $K_5^{(3)}$ with one copy of $L_5$, we decrease the number of edges by $3|E(K_5^{(3)})|-|E(L_5)|=7$. Furthermore, since $L_5$ is aggressively saturated, we can do this replacement repeatedly, giving us constructions on $2n-7k$ vertices for a large range of values $k$. To achieve edge counts $m$ where $2n-m$ is not divisible by $7$, we will need to adjust our construction to match the remainder of $2n-m$ when divided by $7$. As such, we will need $6$ new small constructions, one for each possible non-zero value of $(2n-m \mod 7)$. In order to perfectly replace copies of $K_5^{(3)}$, we will need each of our constructions to have $5k$ vertices from some $k\in\mathbb{N}$. See Table~\ref{table of cases for L=5 upper range} for the specific graphs used.
	
	Before we present each construction, the following table summarizes which case each will be used for. We note that in addition to the $3$-graphs presented below, we also need $S_5$ from Construction~\ref{construction S}.
	
	\begin{center}
		\begin{table}
			\begin{tabular}{|c|c|c|c|c|c|}
				\hline
				$(2n-m \mod 7)$&$3$-Graph $G$&Reference&$|E(G)|$&$|V(G)|$&\begin{minipage}{3.4cm}edge diff between $G$ and $\frac{|V(G)|}{5}K_5^{(3)}$\end{minipage}\\
				\hline
				$1$&$S_5\sqcup K_4^{(3)}$&Construction~\ref{construction S}&$12$&$10$&$-8$\\
				$2$&$R$&Construction~\ref{construction R}&21&$15$&$-9$\\
				$3$&$2B$&Construction~\ref{construction half lantern}&$30$&$20$&$-10$\\
				$4$&$Q$&Construction~\ref{construction Q}&$29$&$20$&$-11$\\
				$5$&$B$&Construction~\ref{construction half lantern}&$15$&$10$&$-5$\\
				$6$&$D$&Construction~\ref{construction D}&$14$&$10$&$-6$\\
				\hline
			\end{tabular}
			\caption{The Constructions necessary for the proof of Theorem~\ref{theorem upper range for L=5}. Diagrams of the graphs featured here can be found in Figures~\ref{figure LandS} and~\ref{figure constructions for upper range l=5}.}\label{table of cases for L=5 upper range}
		\end{table}
	\end{center}

\begin{figure}
		\begin{center}
			\begin{subfigure}[t]{1\textwidth}
				\begin{center}
					\begin{tikzpicture}
						
						\draw[thick] (-0.4,0)--(0,.6928)--(.4,0)--(-0.4,0);
						\draw[thick] (2.6,0)--(3,.6928)--(3.4,0)--(2.6,0);
						\draw[thick] (5.6,0)--(6,.6928)--(6.4,0)--(5.6,0);
						
						\draw[dashed] (-1.3,-0.5)--(1.3,-0.5)--(1.3,1.3)--(-1.3,1.3)--(-1.3,-0.5);
						\draw[dashed] (1.7,-0.5)--(4.3,-0.5)--(4.3,1.3)--(1.7,1.3)--(1.7,-0.5);
						\draw[dashed] (4.7,-0.5)--(7.3,-0.5)--(7.3,1.3)--(4.7,1.3)--(4.7,-0.5);
						
						\draw[dashed] (-1.3,-0.5)--(0,-1.133)--(1.3,-0.5);
						\draw[dashed] (-1.3,1.3)--(0,1.8)--(1.3,1.3);
						
						\draw[dashed] (1.7,-0.5)--(3,-1.133)--(4.3,-0.5);
						\draw[dashed] (1.7,1.3)--(3,1.8)--(4.3,1.3);
						
						\draw[dashed] (4.7,-0.5)--(6,-1.133)--(7.3,-0.5);
						\draw[dashed] (4.7,1.3)--(6,1.8)--(7.3,1.3);

						\node (v1) at (0,1.8) {};
						\node (v2) at (3,1.8) {};
						\node (v3) at (6,1.8) {};
						\node (x11) at (-0.4,0) {};
						\node (x12) at (0,.6928) {};
						\node (x13) at (0.4,0) {};
						\node (x21) at (2.6,0) {};
						\node (x22) at (3,.6928) {};
						\node (x23) at (3.4,0) {};
						\node (x31) at (5.6,0) {};
						\node (x32) at (6,.6928) {};
						\node (x33) at (6.4,0) {};
						\node (u1) at (0,-1.133) {};
						\node (u2) at (3,-1.133) {};
						\node (u3) at (6,-1.133) {};
						\node (a) at (1.5,3.8) {};
						\node (b1) at (0,2.8) {};
						\node (b2) at (3,2.8) {};
						\node (c1) at (6,3.8) {};
						\node (c2) at (6,2.8) {};
						\node (center1) at (0,.2309) {};
						\node (center2) at (3,.2309) {};
						\node (center3) at (6,.2309) {};
						
						\draw (center1) circle (0.6) node {};
						\draw (center2) circle (0.6) node {};
						\draw (center3) circle (0.6) node {};
						
						\fill (v1) circle (0.1) node [below] {$v_1$};
						\fill (v2) circle (0.1) node [below] {$v_2$};
						\fill (v3) circle (0.1) node [below] {$v_3$};
						\fill (x11) circle (0.1) node [below left] {$x_{1,1}$};
						\fill (x12) circle (0.1) node [above] {$x_{1,2}$};
						\fill (x13) circle (0.1) node [below right] {$x_{1,3}$};
						\fill (x21) circle (0.1) node [below left] {$x_{2,1}$};
						\fill (x22) circle (0.1) node [above] {$x_{2,2}$};
						\fill (x23) circle (0.1) node [below right] {$x_{2,3}$};
						\fill (x31) circle (0.1) node [below left] {$x_{3,1}$};
						\fill (x32) circle (0.1) node [above] {$x_{3,2}$};
						\fill (x33) circle (0.1) node [below right] {$x_{3,3}$};
						\fill (u1) circle (0.1) node [below] {$u_1$};
						\fill (u2) circle (0.1) node [below] {$u_2$};
						\fill (u3) circle (0.1) node [below] {$u_3$};
						\fill (a) circle (0.1) node [above] {$a$};
						\fill (b1) circle (0.1) node [below] {$b_1$};
						\fill (b2) circle (0.1) node [below] {$b_2$};
						\fill (c1) circle (0.1) node [below] {$c_1$};
						\fill (c2) circle (0.1) node [below] {$c_2$};
						
						\draw ($(a)+(0,0.5)$)
						to[out=0,in=135] ($(b2)+(0,0.6)$)
						to[out=315,in=90] ($(b2)+(0.3,0)$)
						to[out=270,in=90] ($(v2)+(0.3,0)$)
						to[out=270,in=0] ($(v2)-(0,0.5)$)
						to[out=180,in=270] ($(v2)-(0.3,0)$)
						to[out=90,in=270] ($(b2)-(0.3,0.2)$)
						to[out=90,in=345] ($(a)-(0,0.3)$)
						to[out=165,in=270] ($(a)-(0.3,0)$)
						to[out=90,in=180] ($(a)+(0,0.5)$);
						
						\draw ($(a)+(0,0.5)$)
						to[out=180,in=45] ($(b1)+(0,0.6)$)
						to[out=225,in=90] ($(b1)-(0.3,0)$)
						to[out=270,in=90] ($(v1)-(0.3,0)$)
						to[out=270,in=180] ($(v1)-(0,0.5)$)
						to[out=0,in=270] ($(v1)+(0.3,0)$)
						to[out=90,in=270] ($(b1)-(-0.3,0.2)$)
						to[out=90,in=195] ($(a)-(0,0.3)$)
						to[out=15,in=270] ($(a)+(0.3,0)$)
						to[out=90,in=0] ($(a)+(0,0.5)$);
						
						\draw ($(c1)+(0,0.3)$)
						to[out=0,in=90] ($(c1)+(0.4,0)$)
						to[out=270,in=90] ($(v3)+(0.4,0)$)
						to[out=270,in=0] ($(v3)-(0,0.5)$)
						to[out=180,in=270] ($(v3)-(0.4,0)$)
						to[out=90,in=270] ($(c1)-(0.4,0)$)
						to[out=90,in=180] ($(c1)+(0,0.3)$);
						
						\draw ($(u1)-(0,0.5)$)
						to[out=0,in=180] ($(u3)-(0,0.5)$)
						to[out=0,in=270] ($(u3)+(0.4,0)$)
						to[out=90,in=0] ($(u3)+(0,0.3)$)
						to[out=180,in=0] ($(u1)+(0,0.3)$)
						to[out=180,in=90] ($(u1)-(0.4,0)$)
						to[out=270,in=180] ($(u1)-(0,0.5)$);
						
						\node at (8.25,0.8) {\LARGE{$\cup$}};

						\begin{scope}[xshift=9cm,yshift=-2cm]
							\node (acopy) at (1.7,4.5) {};
							\node (b1copy) at (0.4,2.8) {};
							\node (b2copy) at (3,2.8) {};
							\node (c1copy) at (5.2,4.5) {};
							\node (c2copy) at (5.2,2.8) {};
							\fill (acopy) circle (0.1) node [above] {$a$};
							\fill (b1copy) circle (0.1) node [below] {$b_1$};
							\fill (b2copy) circle (0.1) node [below] {$b_2$};
							\fill (c1copy) circle (0.1) node [below] {$c_1$};
							\fill (c2copy) circle (0.1) node [below] {$c_2$};
							
							\draw ($(c1copy)+(0,0.3)$)
							to[out=0,in=90] ($(c1copy)+(0.4,0)$)
							to[out=270,in=90] ($(c2copy)+(0.4,0)$)
							to[out=270,in=0] ($(c2copy)-(0,0.5)$)
							to[out=180,in=0] ($(b2copy)-(0,0.5)$)
							to[out=180,in=270] ($(b2copy)-(0.4,0)$)
							to[out=90,in=180] ($(b2copy)+(0,0.3)$)
							to[out=0,in=180] ($(c2copy)+(-0.6,0.3)$)
							to[out=0,in=270] ($(c1copy)-(0.4,0)$)
							to[out=90,in=180] ($(c1copy)+(0,0.3)$);
							
							\draw ($(c1copy)+(0,0.45)$)
							to[out=0,in=90] ($(c1copy)+(0.7,0)$)
							to[out=270,in=90] ($(c2copy)+(0.7,0)$)
							to[out=270,in=270] ($(b1copy)-(0.4,0)$)
							to[out=90,in=180] ($(b1copy)+(0,0.3)$)
							to[out=0,in=90] ($(b1copy)+(0.4,0)$)
							to[out=270,in=270] ($(c2copy)-(0.7,0)$)
							to[out=90,in=270] ($(c1copy)-(0.7,0)$)
							to[out=90,in=180] ($(c1copy)+(0,0.45)$);
							
							\draw ($(c1copy)+(0,0.65)$)
							to[out=0,in=90] ($(c1copy)+(0.55,0)$)
							to[out=270,in=0] ($(c1copy)-(0,0.5)$)
							to[out=180,in=0] ($(acopy)-(0,0.5)$)
							to[out=180,in=90] ($(b1copy)+(0.55,0)$)
							to[out=270,in=0] ($(b1copy)-(0,0.55)$)
							to[out=180,in=270] ($(b1copy)-(0.55,0)$)
							to[out=90,in=180] ($(acopy)+(0,0.65)$)
							to[out=0,in=180] ($(c1copy)+(0,0.65)$);
							
							\draw ($(acopy)+(0,0.5)$)
							to[out=0,in=90] ($(acopy)+(0.4,0)$)
							to[out=270,in=180] ($(b2copy)+(0,0.45)$)
							to[out=0,in=180] ($(c2copy)+(0,0.45)$)
							to[out=0,in=90] ($(c2copy)+(0.55,0)$)
							to[out=270,in=0] ($(c2copy)-(0,0.65)$)
							to[out=180,in=0] ($(b2copy)-(0,0.65)$)
							to[out=180,in=270] ($(acopy)-(0.4,0)$)
							to[out=90,in=180] ($(acopy)+(0,0.5)$);
						\end{scope}

					\end{tikzpicture}
				\end{center}
				\caption{The $3$-graph $Q$ from Construction~\ref{construction Q}. For readability, the graph is presented as the union of two sets of edges on the same vertex set.}\label{subfigure Q}
			\end{subfigure}\\
			
			\begin{subfigure}[t]{\textwidth}
				\begin{center}
					\begin{tikzpicture}

						\draw[thick] (1,0)--(0,0)--(0.5,.866)--(1,0);
						\draw[thick] (2,0) -- (2,.866);
						\draw[thick] (5,0)--(6,0)--(5.5,.866)--(5,0);
						\draw[thick] (4,0) -- (4,.866);
						
						\draw[dashed] (-0.6,-0.6)--(2.6,-0.6)--(2.6,1.466)--(-0.6,1.466)--(-0.6,-0.6);
						
						\draw[dashed] (3.4,-0.6)--(6.6,-0.6)--(6.6,1.466)--(3.4,1.466)--(3.4,-0.6);
						
						\draw[dashed] (-0.6,-0.6)--(1,-1.5)--(2.6,-0.6);
						\draw[dashed] (-0.6,1.466)--(1,2.366)--(2.6,1.466);
						\draw[dashed] (3.4,-0.6)--(5,-1.5)--(6.6,-0.6);
						\draw[dashed] (3.4,1.466)--(5,2.366)--(6.6,1.466);
						
						\node (a1) at (0,0) {};
						\node (a3) at (1,0) {};
						\node (a2) at (0.5,.866) {};
						\node (b1) at (2,.866) {};
						\node (b2) at (2,0) {};
						\node (d2) at (4,0) {};
						\node (c1) at (5,0) {};
						\node (c2) at (5.5,.866) {};
						\node (c3) at (6,0) {};
						\node (d1) at (4,.866) {};
						\node (z) at (3,-1) {};
						\node (x1) at (1,2.366) {};
						\node (x2) at (1,-1.5) {};
						\node (y1) at (5,2.366) {};
						\node (y2) at (5,-1.5) {};
						\node (center1) at (0.5,.2887) {};
						\node (center2) at (5.5,.2887) {};
						
						\draw (center1) circle (0.7) node {};
						\draw (center2) circle (0.7) node {};
						
						\fill (a1) circle (0.1) node [below left] {$a_1$};
						\fill (a2) circle (0.1) node [above] {$a_2$};
						\fill (a3) circle (0.1) node [below right] {$a_3$};
						\fill (b1) circle (0.1) node [right] {$b_1$};
						\fill (b2) circle (0.1) node [below] {$b_2$};
						\fill (c1) circle (0.1) node [below left] {$c_1$};
						\fill (c2) circle (0.1) node [above] {$c_2$};
						\fill (c3) circle (0.1) node [below right] {$c_3$};
						\fill (d1) circle (0.1) node [left] {$d_1$};
						\fill (d2) circle (0.1) node [below] {$d_2$};
						\fill (x1) circle (0.1) node [above] {$x_1$};
						\fill (x2) circle (0.1) node [below] {$x_2$};
						\fill (y1) circle (0.1) node [above] {$y_1$};
						\fill (y2) circle (0.1) node [below] {$y_2$};
						\fill (z) circle (0.1) node [below] {$z$};

						\begin{scope}
							\draw ($(b1)+(0,0.5)$)
							to[out=0, in=90] ($(b1)+(0.5,0)$)
							to[out=270, in=135] ($(b2)+(0.5,-0.3)$)
							to[out=315,in=90] ($(z)+(0.4,0)$)
							to[out=270,in=0] ($(z)-(0.2,0.5)$)
							to[out=180,in=270] ($(b2)-(0.3,-0.3)$)
							to[out=90,in=180] ($(b1)+(0,0.5)$);
							
							\draw ($(d1)+(0,0.5)$)
							to[out=180, in=90] ($(d1)+(-0.5,0)$)
							to[out=270, in=45] ($(d2)+(-0.5,-0.3)$)
							to[out=225,in=90] ($(z)+(-0.4,0)$)
							to[out=270,in=180] ($(z)-(-0.2,0.5)$)
							to[out=0,in=270] ($(d2)-(-0.5,-0.3)$)
							to[out=90,in=0] ($(d1)+(0,0.5)$);
							
							\draw ($(b1)+(0,0.4)$)
							to[out=0, in=180] ($(d1)+(0,0.4)$)
							to[out=0, in=90] ($(d1)+(0.5,0)$)
							to[out=270, in=0] ($(d1)-(0,0.4)$)
							to[out=180, in=45]($(d1)-(1.3,0.4)$)
							to[out=225, in=0] ($(b2)-(0,0.6)$)
							to[out=180,in=270] ($(b2)+(-0.4,0.5)$)
							to[out=90,in=180] ($(b1)+(0,0.4)$);
						\end{scope}
						
					\end{tikzpicture}
				\end{center}
				\caption{The $3$-graph $R$ from Construction~\ref{construction R}.}\label{subfigure R}   
			\end{subfigure}\\
			\begin{subfigure}[t]{0.5\textwidth}
				\begin{center}
					\begin{tikzpicture}
						
						\draw[thick] (0.6,0)--(1,.6928)--(1.4,0)--(0.6,0);
						\draw[thick] (3.6,0)--(4,.6928)--(4.4,0)--(3.6,0);
						
						\draw[dashed] (-0.15,-0.5)--(2.15,-0.5)--(2.15,1.3)--(-0.15,1.3)--(-0.15,-0.5);
						\draw[dashed] (2.9,-0.5)--(5.1,-0.5)--(5.1,1.3)--(2.9,1.3)--(2.9,-0.5);
						
						\draw[dashed] (-0.15,1.3)--(0,2.1)--(2.15,1.3);
						\draw[dashed] (-0.15,1.3)--(2,2.1)--(2.15,1.3);
						
						\draw[dashed] (2.9,-0.5)--(4,-1.433)--(5.1,-0.5);
						\draw[dashed] (2.9,1.3)--(4,2.1)--(5.1,1.3);

						\node (v1) at (0,2.1) {};
						\node (v2) at (2,2.1) {};
						\node (a1) at (4,2.1) {};
						\node (x1) at (0.6,0) {};
						\node (x2) at (1,.6928) {};
						\node (x3) at (1.4,0) {};
						\node (y1) at (3.6,0) {};
						\node (y2) at (4,.6928) {};
						\node (y3) at (4.4,0) {};
						\node (a2) at (4,-1.433) {};
						\node (center1) at (1,.2309) {};
						\node (center2) at (4,.2309) {};
						
						\draw (center1) circle (0.6) node {};
						\draw (center2) circle (0.6) node {};
						
						\fill (v1) circle (0.1) node [above] {$v_1$};
						\fill (v2) circle (0.1) node [above] {$v_2$};
						\fill (a1) circle (0.1) node [above] {$a_1$};
						\fill (x1) circle (0.1) node [below left] {$x_{1}$};
						\fill (x2) circle (0.1) node [above] {$x_{2}$};
						\fill (x3) circle (0.1) node [below right] {$x_{3}$};
						\fill (y1) circle (0.1) node [below left] {$y_{1}$};
						\fill (y2) circle (0.1) node [above] {$y_{2}$};
						\fill (y3) circle (0.1) node [below right] {$y_{3}$};
						\fill (a2) circle (0.1) node [below] {$a_2$};
						
						\begin{scope}
							\draw ($(v1)-(0.4,0)$)
							to[out=90,in=180] ($(v1)+(0,0.5)$)
							to[out=0,in=180] ($(a1)+(0,0.5)$)
							to[out=0,in=90] ($(a1)+(0.4,0)$)
							to[out=270,in=0] ($(a1)-(0,0.3)$)
							to[out=180,in=0] ($(v1)-(0,0.3)$)
							to[out=180,in=270] ($(v1)-(0.4,0)$);
						\end{scope}
						
					\end{tikzpicture}
				\end{center}
				\caption{The $3$-graph $B$ from Construction~\ref{construction half lantern}, referred to as the broken lantern.}\label{subfigure B}
			\end{subfigure}~
			\begin{subfigure}[t]{0.5\textwidth}
				\begin{center}
					\begin{tikzpicture}
						
						\draw[dashed] (-0.8,-0.5)--(1.6,-0.5)--(1.6,1.5)--(-0.8,1.5)--(-0.8,-0.5);
						\draw[dashed] (-0.8,1.5)--(0.5,2)--(1.6,1.5);
						
						\draw[dashed] (1.9,-0.5)--(4.3,-0.5)--(4.3,1.5)--(1.9,1.5)--(1.9,-0.5);
						\draw[dashed] (1.9,1.5)--(3,2)--(4.3,1.5);
						
						\draw[thick] (0,0)--(0,1)--(1,1)--(1,0)--(0,0)--(1,1);
						\draw[thick] (1,0)--(0,1);
						
						\draw[thick] (2.5,0)--(2.5,1)--(3.5,1)--(3.5,0)--(2.5,0)--(3.5,1);
						\draw[thick] (2.5,1)--(3.5,0);

						\node (a) at (0.5,2) {};
						\node (x1) at (0,1) {};
						\node (x2) at (1,1) {};
						\node (x3) at (0,0) {};
						\node (x4) at (1,0) {};
						\node (b) at (3,2) {};
						\node (y1) at (2.5,1) {};
						\node (y2) at (3.5,1) {};
						\node (y3) at (2.5,0) {};
						\node (y4) at (3.5,0) {};
						\node (invisible) at (0,-1.3) {};
						
						\fill (a) circle (0.1) node [above] {$a$};
						\fill (x1) circle (0.1) node [left] {$x_1$};
						\fill (x2) circle (0.1) node [right] {$x_2$};
						\fill (x3) circle (0.1) node [left] {$x_3$};
						\fill (x4) circle (0.1) node [right] {$x_4$};
						\fill (b) circle (0.1) node [above] {$b$};
						\fill (y1) circle (0.1) node [left] {$y_1$};
						\fill (y2) circle (0.1) node [right] {$y_2$};
						\fill (y3) circle (0.1) node [left] {$y_3$};
						\fill (y4) circle (0.1) node [right] {$y_4$};
						\fill[white] (invisible) circle (0.1) node {};

						\begin{scope}
							\draw ($(x1)-(0.65,0)$)
							to[out=90,in=180] ($(x1)+(0,0.35)$)
							to[out=0, in=180] ($(y1)+(0,0.35)$)
							to[out=0,in=90] ($(y1)+(0.35,0)$)
							to[out=270,in=0] ($(y1)-(0,0.35)$)
							to[out=180,in=0] ($(x1)-(0,0.35)$)
							to[out=180,in=270] ($(x1)-(0.65,0)$);
							
							\draw ($(x4)-(0.35,0)$)
							to[out=90,in=180] ($(x4)+(0,0.35)$)
							to[out=0,in=180] ($(y4)+(0,0.35)$)
							to[out=0,in=90] ($(y4)+(0.65,0)$)
							to[out=270,in=0] ($(y4)-(0,0.35)$)
							to[out=180,in=0] ($(x4)-(0,0.35)$)
							to[out=180,in=270] ($(x4)-(0.35,0)$);
						\end{scope}
					\end{tikzpicture}
				\end{center}
				\caption{The $3$-graph $D$ from Construction~\ref{construction D}.}\label{subfigure D}
			\end{subfigure}
		\end{center}
		\caption{Diagrams of some of the constructions required for the proof of Theorem~\ref{theorem upper range for L=5}. A dashed rectangle represents that the $2$-graph inside the rectangle is part of the link of any vertex connected to the rectangle via dashed lines. All smooth curves represent $3$-edges in the graph.}\label{figure constructions for upper range l=5}
	\end{figure}
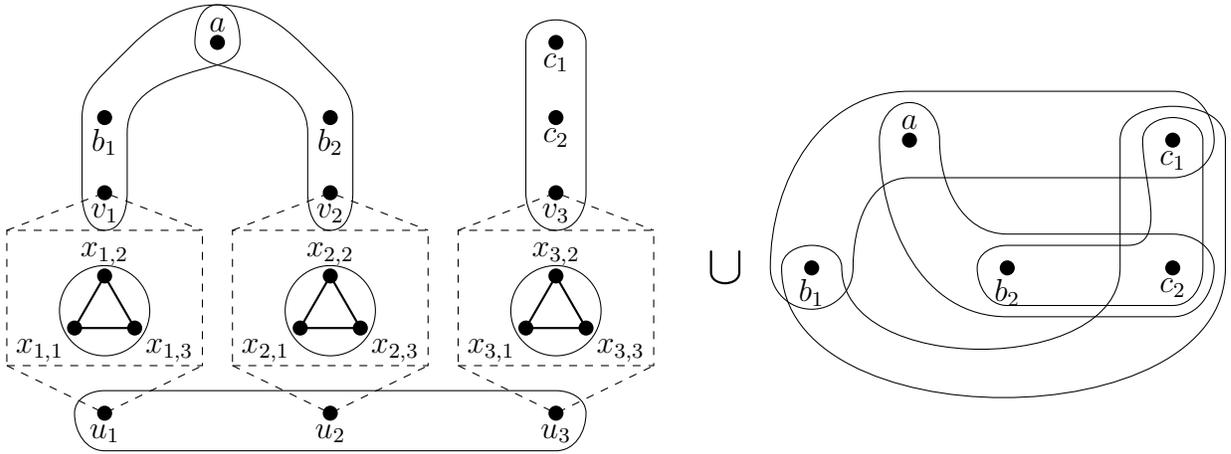
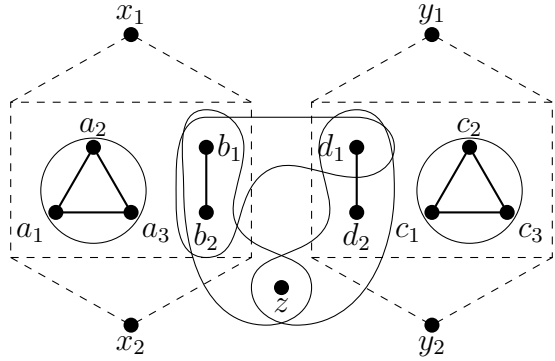
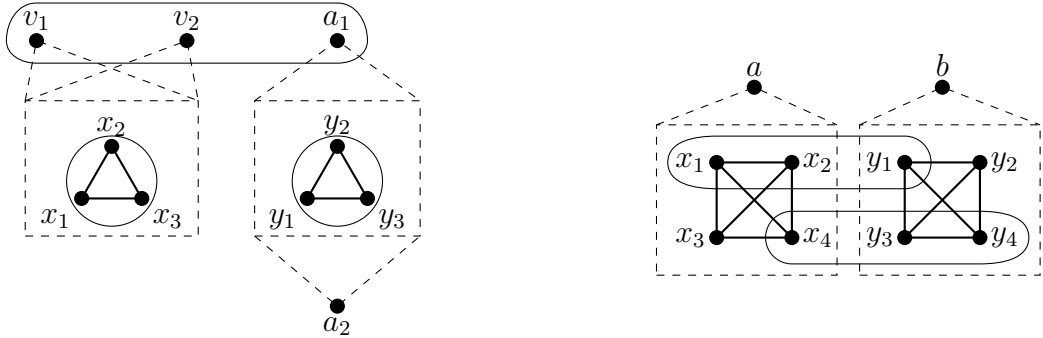

	\begin{construction}\label{construction half lantern}
		Let $B$ denote the \textbf{broken lantern}, i.e. the graph on $10$ vertices with
		\[
		V(B)=\{a_1,a_2,v_1,v_2,x_i,y_i\mid i\in [3]\}
		\]
		and
		\[
		E(B)=\{a_1v_1v_2, x_1x_2x_3,y_1y_2y_3, a_ix_jx_{j'},v_iy_jy_{j'}\mid i\in [2],j,j'\in [3], j\neq j'\}.
		\]
		We note that $|E(B)|=15$. See Figure~\ref{subfigure B} for a diagram of the broken lantern.
	\end{construction}
	
	\begin{construction}\label{construction D}
		Let $D$ denote the graph on $10$ vertices with
		\[
		V(D)=\{a,b,x_i,y_i\mid i\in [4]\}
		\]
		and
		\[
		E(D)=\{ax_ix_j,by_iy_j,x_1x_2y_1,y_3y_4x_4\mid i,j\in [4],i\neq j\}.
		\]
		We note that $|E(D)|=14$. See Figure~\ref{subfigure D} for a diagram of $D$.
	\end{construction}
	
	\begin{construction}\label{construction Q}
		Let $Q$ denote the graph on $20$ vertices formed from the $5$-lantern $L_5$ as follows. Let
		\[
		V(Q)=V(L_5)\sqcup \{a,b_1,b_2,c_1,c_2\}
		\]
		and let
		\[
		E(Q)=(E(L_5)\setminus\{v_1v_2v_3\})\cup\{ab_1v_1,ab_2v_2,ab_1c_1,ab_2c_2,c_1c_2v_3,b_1c_1c_2,b_2c_1c_2\}.
		\]
		We note that $|E(Q)|=29$. See Figure~\ref{subfigure Q} for a diagram of $Q$.
	\end{construction}
	
	\begin{construction}\label{construction R}
		Let $R$ denote the graph on $15$ vertices with
		\[
		V(R)=\{x_i,y_i,z,a_j,b_i,c_j,d_i\mid i\in [2],j\in [3]\},
		\]
		and with edges
		\[
		E(R)=\{x_ia_ja_{j'},x_ib_1b_2,y_ic_jc_{j'},y_id_1d_2,a_1a_2a_3,c_1c_2c_3,b_1b_2z,d_1d_2z,b_1b_2d_1\mid i\in [2],j,j'\in [3], j\neq j'\}.
		\]
		We note that $|E(R)|=21$. See Figure~\ref{subfigure R} for a diagram of $R$.
	\end{construction}
	
	\begin{proposition}\label{proposition D Q R B are all saturated}
		The graphs $D, Q$ and $R$ from Constructions\ref{construction D}~\ref{construction Q} and~\ref{construction R} are all Berge-$K_{1,5}$-saturated. Furthermore, if $B$ is the graph from Construction~\ref{construction half lantern}, then both $B$ and $2B$ are Berge-$K_{1,5}$-saturated.
	\end{proposition}
	
	\begin{proof}
		As each graph is small, it is readily verified that none of these graphs contain a Berge-$K_{1,5}$. We will consider each graph individually to show that adding an edge creates a Berge-$K_{1,5}$.
		
		\textbf{$D$:} We have that $L(a)\cong L(b)\cong K_4$, while $L(y_3)\cong L(x_3)\cong K_{1,3}\sqcup K_2$. Therefore, all these vertices are aggressively saturated of type I. The links of vertices in $\{x_1,x_2,y_1,y_2\}$ consist of a triangle with a pendant edge, where the degree $1$ vertex in these links are either $x_3$ or $y_3$, and therefore the missing edges in this link all include an aggressively saturated vertex of type I, and thus the vertices in $\{x_1,x_2,y_1,y_2\}$ are aggressively saturated of type II. This only leaves $x_4$ and $y_4$, as non-aggressively saturated. Since two vertices trivially induce a clique, Corollary~\ref{corollary if non agg sat induce clique} gives us that $D$ is saturated.
		
		\textbf{$Q$:} We can see that for every vertex in $\{u_i,v_i\mid i\in [3]\}$, the links are isomorphic to $K_3\sqcup K_2$, and therefore all these vertices are aggressively saturated of type I. Also, the links of all the vertices in $\{x_{i,j}\mid i,j\in [3]\}$ are $K_4^-$, where the missing edge in these links contains one of the $u_i$'s, and therefore these vertices are aggressively saturated of type II. We further see that the link of $a$ is $2P_3$, while the links of $c_1$ and $c_2$ are isomorphic to $T_0$ (the unique tree of order $5$ with a degree $3$ vertex), and thus $a$, $c_1$ and $c_2$ are all aggressively saturated of type I. This only leaves $b_1$ and $b_2$ as non-aggressively saturated vertices. Corollary~\ref{corollary if non agg sat induce clique} gives us that $Q$ is saturated.
		
		\textbf{$R$:} We have that the links of all the vertices in $\{x_1,x_2,y_1,y_2\}$ are isomorphic to $K_3\sqcup K_2$, and therefore these vertices are aggressively saturated of type I. Vertices in $\{a_1,a_2,a_3,c_1,c_2,c_3\}$ have links isomorphic to $K_4^-$, where the missing edge in these links involves a vertex in $\{x_1,x_2,y_1,y_2\}$, and therefore these vertices are all aggressively saturated of type II. Vertices $b_1$ and $b_2$ have links $K_{1,4}$, while $d_1$ has link $K_{1,3}\sqcup K_2$, and therefore these three vertices are aggressively saturated of type I. This leaves only $d_2$ and $z$ as non-aggressively saturated, and therefore Corollary~\ref{corollary if non agg sat induce clique} gives us that $R$ is saturated.
		
		\textbf{$2B$:} Let us consider a single copy of $B$ for a moment. We find that $L(a_1)\cong L(v_1)\cong L(v_2)\cong K_3\sqcup K_2$, and therefore these vertices are aggressively saturated of type I in $B$. Furthermore, $L(x_i)\cong L(y_i)\cong K_4^-$ for all $i\in [3]$, where the missing edge from each of these $K_4^-$'s contains one of the vertices in $\{a_1,v_1,v_2\}$, and therefore each of these vertices is aggressively saturated of type II. This leaves $a_2$ as the only non-aggressively saturated vertex in $B$.
		
		Now, in $2B$, all the vertices except the two copies of $a_2$ are aggressively saturated, and trivially two vertices form a clique. By Corollary~\ref{corollary if non agg sat induce clique}, $2B$ is saturated.
		
		\textbf{$B$:} Since $2B$ is saturated, every component of $2B$ must be saturated, and thus $B$ itself is saturated.
	\end{proof}

	\subsection{The precise upper range for \texorpdfstring{$\ell=5$}{L=5}}
	
	We now provide the final result to completely determine the spectrum for Berge-$K_{1,5}$ when $5\mid n$.
	
	\begin{theorem}\label{theorem upper range for L=5}
		For all $n$ divisible by $5$ with $n$ large enough, and for all $m\in [\frac{5}{3}n,2n-5]$, there exists an $n$ vertex Berge-$K_{1,5}$-saturated graph with $m$ edges.
	\end{theorem}
	
	\begin{proof}
		Let $m\in [\frac{5}{3}n,2n-5]$, and let $m^*:=2n-m$. Let $a,b\in\mathbb{Z}_{\geq 0}$ with $0\leq b<7$ be such that $m^*=7a+b$. We now describe a construction $G$ with $n$ vertices and $m$ edges based on $b$.
		
		\textbf{Case 1:} $b=0$. Then let 
		\[
		G:=aL_5\sqcup \left(\frac{n}{5}-3a\right)K^{(3)}_5.
		\]
		We have that $|V(G)|=15a+5\left(\frac{n}{5}-3a\right)=n$, while 
		\[
		|E(G)|=23a+10\left(\frac{n}{5}-3a\right)=2n-7a=m.
		\]
		
		\textbf{Case 2:} $b=1$. Then let 
		\[
		G:=S_5\sqcup K_4^{(3)}\sqcup (a-1)L_5\sqcup \left(\frac{n}{5}-3a-2\right)K^{(3)}_5.
		\]
		We have that $|V(G)|=6+4+15(a-1)+5\left(\frac{n}{5}-3a+1\right)=n$, while
		\[
		|E(G)|=8+4+23(a-1)+10\left(\frac{n}{5}-3a+1\right)=2n-7a-1=m.
		\]
		
		\textbf{Case 3:} $b=2$. Then let 
		\[
		G:=R\sqcup (a-1)L_5\sqcup \left(\frac{n}{5}-3a\right)K^{(3)}_5.
		\]
		We have that $|V(G)|=15+15(a-1)+5\left(\frac{n}{5}-3a\right)=n$, while
		\[
		|E(G)|=21+23(a-1)+10\left(\frac{n}{5}-3a\right)=2n-7a-2=m.
		\]
		
		\textbf{Case 4:} $b=3$. Then let 
		\[
		G:=2B\sqcup (a-1)L_5\sqcup \left(\frac{n}{5}-3a-1\right)K^{(3)}_5.
		\]
		We have that $|V(G)|=10(2)+15(a-1)+5\left(\frac{n}{5}-3a-1\right)=n$, while
		\[
		|E(G)|=15(2)+23(a-1)+10\left(\frac{n}{5}-3a-1\right)=2n-7a-3=m.
		\]
		
		\textbf{Case 5:} $b=4$. Then let 
		\[
		G:=Q\sqcup (a-1)L_5\sqcup \left(\frac{n}{5}-3a-1\right)K^{(3)}_5.
		\]
		We have that $|V(G)|=20+15(a-1)+5\left(\frac{n}{5}-3a-1\right)=n$, while
		\[
		|E(G)|=29+23(a-1)+10\left(\frac{n}{5}-3a-1\right)=2n-7a-4=m.
		\]
		
		\textbf{Case 6:} $b=5$. Then let 
		\[
		G:=B\sqcup aL_5\sqcup \left(\frac{n}{5}-3a-2\right)K^{(3)}_5.
		\]
		We have that $|V(G)|=10+15a+5\left(\frac{n}{5}-3a-2\right)=n$, while
		\[
		|E(G)|=15+23a+10\left(\frac{n}{5}-3a-2\right)=2n-7a-5=m.
		\]
		
		\textbf{Case 7:} $b=6$. Then let 
		\[
		G:=D\sqcup aL_5\sqcup \left(\frac{n}{5}-3a-2\right)K^{(3)}_5.
		\]
		We have that $|V(G)|=10+15a+5\left(\frac{n}{5}-3a-2\right)=n$, while
		\[
		|E(G)|=14+23a+10\left(\frac{n}{5}-3a-2\right)=2n-7a-6=m.
		\]
		
		We quickly verify that the above quantities are all non-negative. Since $m\geq \frac{5}{3}n$, we have $m^*\leq \frac{n}{3}$, and so $a\leq \frac{n}{21}+1$. In particular, the quantity $\frac{n}{5}-3a+i\geq 0$ for all $i\geq -2$ and $n$ large enough. In addition, as long as $0\leq b\leq 4$, we have that $a-1\geq 0$ since $m^*\geq 5$.
		
		Now, we note that in all cases, $G$ is Berge-$K_{1,5}$-saturated. Indeed, from Lemma~\ref{lemma L is aggressively saturated},~\ref{lemma S is aggressively saturated} and Observation~\ref{observation K is aggressively saturated}, we have that $L_5$, $S_5$ and $K^{(3)}_5$ are all aggressively saturated. Then by Observation~\ref{observation union of aggresively saturated graphs}, disjoint unions of these are also aggressively saturated. Finally, by Proposition~\ref{proposition D Q R B are all saturated}, we have that the graphs $B$, $2B$, $D$, $Q$, and $R$ are all saturated. Then $G$ is always expressible as a disjoint union of an aggressively saturated graph with a saturated graph, and therefore $G$ is saturated.
	\end{proof}
	
	\section{Concluding remarks}
	
	In this work, we were able to determine exactly the $3$-uniform saturation spectrum of Berge-$K_{1,5}$ when $5\mid n$. For Berge-$K_{1,\ell}$ with $\ell>5$, we do not completely determine the upper range of the spectrum. It would be interesting to explore this range, however some new ideas would be needed. In particular, the reason Berge-$K_{1,5}$ is tractable is due to the fact that the number of edges in the $5$-lantern $L_5$ is only $7$ fewer than the number of edges in $3$ copies of $K^{(3)}_{5}$; this allows us to consider $7$ cases in the proof of Theorem~\ref{theorem upper range for L=5}. However, as $\ell$ increases, the difference in the number of edges between $L_{\ell}$ and $3K^{(3)}_\ell$ increases, and this casework becomes unfeasible. While the $\ell=6$ case may be resolvable with considerable effort, it is unclear how to determine the full upper range for all $\ell\geq 6$ without the development of significant new tools.
	
	Another natural direction for future work is to extend these results to hypergraphs of uniformity $k>3$. Many of the ideas here apply equally well in higher uniformities, but new issues arise. As such, substantially more work would be needed to extend our results even to $4$-graphs. In particular, many of our results rely on the analysis of links of vertices, which becomes significantly more difficult when the link graphs are not $2$-graphs.
	
	Finally, it would be interesting to consider the saturation spectrum for other hypergraphs beyond the stars considered in this manuscript. A significant difficulty here is that there are very few choices of hypergraphs for which both the extremal number and saturation number are known exactly. One possible example might be to consider $3$-uniform Berge-$K_3$. In~\cite{EGGMS}, it was shown that $\sat_3(n,\text{Berge-}K_3)=\left\lceil \frac{n-1}{2}\right\rceil$, while in~\cite{Gy}, it was shown that $\ex(n,\text{Berge-}K_3)=\left\lfloor\frac{n^2}{8}\right\rfloor$. It is our suspicion that, as in the case of the star, the saturation spectrum contains all but finitely many values near the top of the range.


\begin{thebibliography}{10}
		\label{bibliography}
		
		
		\bibitem{AE} B. Austhoff, S. English. \emph{Nearly-Regular Hypergraphs and Saturation of Berge Stars}. Electronic J. Combin. {\bf{26}}(4), (2019). 
		
		\bibitem{Amin} K. Amin, J. Faudree, R. J. Gould, E. Sidorowicz. \emph{On the non-$(p-1)$-partite $K_p$-free graphs}. Discuss. Math. Graph Theory, 33(1), pp. 9--23 (2013).
		
		\bibitem{SS2} P. Balister, A. Dogan. \emph{On the edge spectrum of saturated graphs for paths and stars}. J. Graph Theory 84(9), pp. 364--385 (2018).
		
		\bibitem{Barefoot} C. Barefoot, K. Casey, D. Fisher, K. Fraughnaugh, H. Harary. \emph{Size in maximal triangle- free graphs and minimal graphs of diameter 2}. Discrete Math. {\bf{138}}, pp. 93--99 (1995).
		
		\bibitem{RecentEG} S. Behrens, C. Erbes, M. Ferrara, S. Hartke, B. Reiniger, H. Spinoza, C. Tomlinson. \emph{New results on degree sequences of uniform hypergraphs}. Electron. J. Combin. {\bf{20}}(4), pp. 14--18 (2013).
		
		\bibitem{Berge} C. Berge. \emph{Graphs and hypergraphs}. Elsevier Science Ltd. (1985).
		
		\bibitem{LinearEG} N. Bhave, B. Bam, C. Deshpande. \emph{A characterization of degree sequences of linear hypergraphs}. J. Indian Math. Soc. (N.S.) {\bf{75}} no. 1--4, pp. 151--160 (2009).
		
		\bibitem{Bollobas} B. Bollob\'as. \emph{On generalized graphs}. Acta Math. Acad. Sci. Hungar. {\bf{16}}, pp. 447--452 (1965).
		
		\bibitem{Configs} M. Bras-Amor\'os, K. Stokes. \emph{The semigroup of combinatorial configurations}. Semigroup Forum {\bf{84}}(1), pp. 91--96 (2012).
		
		\bibitem{d} D. de~Caen, \emph{Extension of a theorem of Moon and Moser on complete subgraphs}, Ars Combin. {\bf 16}, pp.  5--10 (1983).
		
		\bibitem{EGGMS} S.~English, P.~Gordon, N.~Graber, A.~Methuku, E.~Sullivan. \emph{Saturation of Berge hypergraphs}, Discrete Math. {\bf 342}(6), pp. 1738--1761 (2019).
		
		\bibitem{EKNS} S.~English, J.~Kritschgau, M.~Nahvi, E.~Sprangel. \emph{Saturation numbers for Berge cliques}, European J. Combin. {\bf 118}, 103911 (2024).
		
		\bibitem{EG} P. Erd\H{o}s, T. Gallai. \emph{On maximal paths and circuits of graphs}. Acta Math. Acad. Sci. Hungary {\bf{10}}, pp. 337--356 (1959).
		
		\bibitem{EHM} P.~Erd\H{o}s, A.~Hajnal, J.~Moon. \emph{A problem in graph theory}. Amer. Math. Monthly {\bf{71}}, pp. 1107--1110 (1964).
		
		\bibitem{ES} P.~Erd\H{o}s, A.~Stone. \emph{On the structure of linear graphs}. Bull. Amer. Math. Soc. {\bf{52}}, pp. 1087--1091 (1946).
		
		\bibitem{FL} Z. F\"uredi and R. Luo, \emph{Induced Tur\'an problems and traces of hypergraphs}, European J. Combin. {\bf 111}, 103692 (2023).
		
		\bibitem{SS3} J. Gao, X. Hou, Y. Ma. \emph{The edge spectrum of $K_4^−$-saturated graphs}. Graphs Combin. {\bf{34}}(6), pp. 1723-1727 (2018).
		
		\bibitem{G} D. Gerbner, \emph{The Tur\'an number of Berge book hypergraphs}, SIAM J. Discrete Math. {\bf 38}(4), pp. 2896--2912 (2024).
		
		\bibitem{GMP} D. Gerbner, A. Methuku, and C. Palmer. \emph{General lemmas for Berge-Tur\'an hypergraph problems}. European J. Combin. {\bf{86}},  103082 (2020).
		
		\bibitem{GP} D. Gerbner, C. Palmer. \emph{Extremal results for Berge hypergraphs}. SIAM J. Discrete Math. {\bf{31}}(4), pp. 2314--2327 (2017).
		
		\bibitem{GPTV} D.~Gerbner, B.~Patk\'os, Zs.~Tuza, and M.~Vizer. \emph{On saturation of Berge hypergraphs}, European J. Combin. {\bf 102}, 103477 (2022).
		
		\bibitem{SS4} R. J. Gould, A. K\"undgen, M. Kang. \emph{On the saturation spectrum of odd cycles}. J. Graph Theory {\bf{106}}(2), pp. 213--224 (2024).
		
		\bibitem{SS5} R. J. Gould, W. Tang, E. Wei, C. Q. Zhang. \emph{Edge spectrum of saturation numbers for small paths}. Discrete Math. {\bf{312}}, pp. 2682--2689 (2012).
		
		\bibitem{Gy} E. Gy\H ori. \emph{Triangle-free hypergraphs}, Combin. Probab. Comput. {\bf 15}(1-2), pp. 185--191 (2006).
		
		\bibitem{GS} E. Gy\H ori and N. Salia. \emph{Linear three-uniform hypergraphs with no Berge path of given length}, J. Combin. Theory Ser. B {\bf 171}, pp. 36--48 (2025).
		
		\bibitem{JLR} S. Janson, T. \L uczak and A. Ruci\'nski. {\it Random graphs}, Wiley-Interscience Series in Discrete Mathematics and Optimization, Wiley-Interscience, New York (2000).
		
		\bibitem{KT} L.~K\'aszonyi, Z.~Tuza. \emph{Saturated graphs with minimal number of edges}.
		Journal of Graph Theory {\bf{10}}, pp. 203--210 (1986).
		
		\bibitem{Mantel} W. Mantel. \emph{Problem 28}. {\it Wiskundige Opgaven} {\bf{10}}, pp.60--61 (1907). 
		
		\bibitem{Pikhurko} O. Pikhurko. \emph{The minimum size of saturated hypergraphs}. Combin. Probab. Comput. {\bf{8}}(5), pp. 483--492 (1995).
		
		\bibitem{Turan} P. Tur\'an. {\emph{On an extremal problem in graph theory.}} Mat. Fiz. Lap. (in Hungarian) {\bf{48}}, pp. 436--452 (1941).
		
	\end{thebibliography}
\end{document}